%%%%%%%%%%%%%%%%%%%%%%%%%%%%%%%%%%%%%%%%%%%%%%%%%%%%%%%%%%
%%%%%%%%%%%%%%%%%%%%%%%%%%%%%%%%%%%%%%%%%%%%%%%%%%%%%%%%%%
%%
%%     This is the AMS-LaTeX file:
%%
%%     CGMR5
%%     controllo sul modello del tumore
%%     
%%
%%%%%%%%%%%%%%%%%%%%%%%%%%%%%%%%%%%%%%%%%%%%%%%%%%%%%%%%%%

\def\LaTeX{%
  \let\Begin\begin
  \let\End\end
  \let\salta\relax
  \let\finqui\relax
  \let\futuro\relax}

\def\UK{\def\our{our}\let\sz s}
\def\USA{\def\our{or}\let\sz z}

\UK
%\USA

%%%%%%%%%%%%%%%%%%%%%%%%%%%%%%%%%

% scegliere fra \TeX e \LaTeX  e fra  \UK oppure \USA

%\TeX
\LaTeX

%\UK
\USA

%%%%%%%%%%%%%%%%%%%%%%%%%%%%%%%%%
%% page layout
%%%%%%%%%%%%%%%%%%%%%%%%%%%%%%%%%

\salta

\documentclass[twoside,12pt]{article}
\setlength{\textheight}{24cm}
\setlength{\textwidth}{16cm}
\setlength{\oddsidemargin}{2mm}
\setlength{\evensidemargin}{2mm}
\setlength{\topmargin}{-15mm}
\parskip2mm

%%%%%%%%%%%%%%%%%%%%%%%%%%%%%%%%%
%% packages
%%%%%%%%%%%%%%%%%%%%%%%%%%%%%%%%%

%\usepackage{color}
\usepackage[usenames,dvipsnames]{color}
\usepackage{amsmath}
\usepackage{amsthm}
\usepackage{amssymb, bbm}
\usepackage[mathcal]{euscript}
\usepackage{cite}

\usepackage{hyperref}
%
%		COLORS FOR CORRECTIONS
%
% do the same, please (i.e., don't use the standard {\color{red} text} or similar): 
% just choose the color you prefer in \def\yourname

% example of use:  \juerg{I want this to become blue}

\definecolor{viola}{rgb}{0.3,0,0.7}
\definecolor{ciclamino}{rgb}{0.5,0,0.5}
\definecolor{rosso}{rgb}{0.85,0,0}

\def\revis #1{{\color{cyan}#1}}

\def\pier #1{#1}
\def\gagi #1{#1}
\def\revis #1{#1}

%%%%%%%%%%%%%%%%%%%%%%%%%%%%%%%%%
%% bibliographystyle
%%%%%%%%%%%%%%%%%%%%%%%%%%%%%%%%%

\bibliographystyle{plain}

%%%%%%%%%%%%%%%%%%%%%%%%%%%%%%%%%
%% environments
%%%%%%%%%%%%%%%%%%%%%%%%%%%%%%%%%

%

\finqui

\def\Beq{\Begin{equation}}
\def\Eeq{\End{equation}}

\def\Bthm{\Begin{theorem}}
\def\Ethm{\End{theorem}}
\def\Blem{\Begin{lemma}}
\def\Elem{\End{lemma}}
\def\Bprop{\Begin{proposition}}
\def\Eprop{\End{proposition}}

\def\Brem{\Begin{remark}\rm}
\def\Erem{\End{remark}}

\def\Bdim{\Begin{proof}}
\def\Edim{\End{proof}}
\def\Bcenter{\Begin{center}}
\def\Ecenter{\End{center}}
\let\non\nonumber

%%%%%%%%%%%%%%%%%%%%%%%%%%%%%%%%%
%% macros
%%%%%%%%%%%%%%%%%%%%%%%%%%%%%%%%%

% macro salvate

% sottosezioni non numerate

\def\step #1 \par{\medskip\noindent{\bf #1.}\quad}

% abbreviazioni di parole

\def\Lip{Lip\-schitz}
\def\Holder{H\"older}
\def\Frechet{Fr\'echet}
\def\Carath{Carath\'eodory}
\def\aand{\quad\hbox{and}\quad}

\def\lhs{left-hand side}
\def\rhs{right-hand side}

% bold, cal e mathop

\def\multibold #1{\def\arg{#1}%
  \ifx\arg\pto \let\next\relax
  \else
  \def\next{\expandafter
    \def\csname #1#1#1\endcsname{{\boldsymbol #1}}%
    \multibold}%
  \fi \next}

\def\pto{.}

\def\multical #1{\def\arg{#1}%
  \ifx\arg\pto \let\next\relax
  \else
  \def\next{\expandafter
    \def\csname cal#1\endcsname{{\cal #1}}%
    \multical}%
  \fi \next}

\def\multigrass #1{\def\arg{#1}%
  \ifx\arg\pto \let\next\relax
  \else
  \def\next{\expandafter
    \def\csname gr#1\endcsname{{\mathbb #1}}%
    \multigrass}%
  \fi \next}

% operatori

\def\multimathop #1 {\def\arg{#1}%
  \ifx\arg\pto \let\next\relax
  \else
  \def\next{\expandafter
    \def\csname #1\endcsname{\mathop{\rm #1}\nolimits}%
    \multimathop}%
  \fi \next}

\multibold
qwertyuiopasdfghjklzxcvbnmQWERTYUIOPASDFGHJKLZXCVBNM.

\multical
QWERTYUIOPASDFGHJKLZXCVBNM.

\multigrass
QWERTYUIOPASDFGHJKLZXCVBNM.

\multimathop
diag dist div dom mean meas sign supp .

\def\Span{\mathop{\rm span}\nolimits}

% accorpamenti di formule citate:
% uso  \accorpa {prima}{seconda}
%      \Accorpa\cs prima seconda (con il comodo blank anche dopo)
% NB: \Accorpa definisce \cs come l'accorpamento delle due citazioni
% e scrive sul file.log

\def\accorpa #1#2{\eqref{#1}--\eqref{#2}}
\def\Accorpa #1#2 #3 {\gdef #1{\eqref{#2}--\eqref{#3}}%
  \wlog{}\wlog{\string #1 -> #2 - #3}\wlog{}}

% macro comode

\def\separa{\noalign{\allowbreak}}

\def\somma #1#2#3{\sum_{#1=#2}^{#3}}

\def\graffe #1{\mathopen\{#1\mathclose\}}

\def\<#1>{\mathopen\langle #1\mathclose\rangle}
\def\norma #1{\mathopen \| #1\mathclose \|}

\def\[#1]{\mathopen\langle\!\langle #1\mathclose\rangle\!\rangle}

\def\ioT {\int_0^T}
\def\intQt{\int_{Q_t}}
\def\intQ{\int_Q}
\def\iO{\int_\Omega}

\def\dt{\partial_t}
\def\dn{\partial_\nnn}
\def\ddt{\frac d{dt}}

\def\cpto{\,\cdot\,}

\def\checkmmode #1{\relax\ifmmode\hbox{#1}\else{#1}\fi}

\def\aeQ{\checkmmode{a.e.\ in~$Q$}}

\def\aet{\checkmmode{a.e.\ in~$(0,T)$}}

\def\aaQ{\checkmmode{for a.a.~$(x,t)\in Q$}}

\def\aat{\checkmmode{for a.a.~$t\in(0,T)$}}

% insiemi numerici

\def\erre{{\mathbb{R}}}

% spazi di funzioni a valori vettoriali su [0,T], [0,t], [0,s], [0,+\infty), [\delta,T]

% Come ricordare: in generale i simboli L H W  C da soli per gli spazi su (0,T)
% gli stessi raddoppiati per (0,+\infty)
% aggiunta di t o s al simbolo per (0,t) e (0,s)
% aggiunta di d al simbolo semplice o doppio per intervalli (\delta,T) e (\delta,+\infty)
% il simbolo C e i suoi derivati mettono le quadre anziche' le tonde

% Esempi   \L2V   \L\infty\Vp   \W{1,1}H   \C0H   \LL2V   \CC0\Vp   \Ld2V  \CCdH

\def\genspazio #1#2#3#4#5{#1^{#2}(#5,#4;#3)}
\def\spazio #1#2#3{\genspazio {#1}{#2}{#3}T0}

\def\L {\spazio L}
\def\H {\spazio H}

\def\C #1#2{C^{#1}([0,T];#2)}

% spazi di funzioni su \Omega, \Gamma, Q e \Sigma

\def\Lx #1{L^{#1}(\Omega)}
\def\Hx #1{H^{#1}(\Omega)}

\def\LQ #1{L^{#1}(Q)}

\def\Ldue{\Lx 2}
\def\Linfty{\Lx\infty}

\def\Huno{\Hx 1}
\def\Hdue{\Hx 2}

% lettere greche

\let\badphi\phi
\let\phi\varphi

\let\theta\vartheta

\let\lam\lambda
\let\Lam\Lambda
\let\ka\kappa
\def\Kappa{\hat\ka}

\let\TeXchi\chi                         % new \chi, exactly on the baseline
\newbox\chibox
\setbox0 \hbox{\mathsurround0pt $\TeXchi$}
\setbox\chibox \hbox{\raise\dp0 \box 0 }
\def\chi{\copy\chibox}

% quadratino di fine dimostrazione

% abbreviazioni specifiche del lavoro

\let\emb\hookrightarrow

\def\cdelta{c_\delta}
\def\OmegaC{\Omega_C}
\def\QC{Q_C}
\def\intQC{\int_{Q_C}}
\def\iOC{\int_{\OmegaC}}

\def\phie{\badphi_e}
\def\phir{\badphi_r}
\def\betai{\beta_i}
\def\betae{\beta_e}
\def\betamax{\beta^*}
\def\gammamax{\gamma^*}
\def\kas{{\kappa_s}}
\def\kae{{\kappa_e}}
\def\kai{{\kappa_i}}
\def\kar{{\kappa_r}}
\def\kamin{\kappa_*}
\def\kamax{\kappa^*}
\def\umin{u_{\rm min}}
\def\umax{u_{\rm max}}

\def\sz{s_0}
\def\ez{e_0}
\def\iz{i_0}
\def\rz{r_0}

 % esempi di uso: \y s1 per s^{(1)}, \y i2 per i^{(2)}
 % esempi di uso: \sol1, \sol2, \sol j

\def\kasstar{\ka_s^*} % controllo ottimo
\def\kaestar{\ka_e^*}
\def\kaistar{\ka_i^*}
\def\karstar{\ka_r^*}
\def\sstar{s^*}  % stato corrispondente
\def\estar{e^*}
\def\istar{i^*}
\def\rstar{r^*}
\def\nstar{n^*}

\def\kaslam{\ka_s^\lam} % controlli incrementati
\def\kaelam{\ka_e^\lam}
\def\kailam{\ka_i^\lam}
\def\karlam{\ka_r^\lam}
\def\slam{s^\lam}
\def\elam{e^\lam}
\def\ilam{i^\lam}
\def\rlam{r^\lam}
\def\nlam{n^\lam}

\def\xilam{\xi_\lam} % rapporti incrementali
\def\etalam{\eta_\lam}
\def\iotalam{\iota_\lam}
\def\rholam{\rho_\lam}

  % variazioni

\def\soluz{(s,e,i,r)}
\def\controllam{(\kaslam,\kaelam,\kailam,\karlam)}
\def\soluzl{(\xi,\eta,\iota,\rho)}
\def\soluzstar{(\sstar,\estar,\istar,\rstar)} % stato ottimo
\def\soluza{(p,q,w,z)} % soluzione aggiunto
\def\control{(\kas,\kae,\kai,\kar)} % controllo generico
\def\controlstar{(\kasstar,\kaestar,\kaistar,\karstar)} % controllo ottimo
\def\soluzlam{(\slam,\elam,\ilam,\rlam)} % stato corrispondente all variazione del controllo

\def\Uad{\calU_{ad}}

\def\Vp{V^*}

\def\normaV #1{\norma{#1}_V}

\let\hat\widehat

%%%%%%%%%%%%%%%%%%%%%%%%%%%%%%
\Begin{document}
%%%%%%%%%%%%%%%%%%%%%%%%%%%%%%%%%

%%%%%%%%%%%%%%%%%%%%%%%%%%%%%%%%%
%% front page
%%%%%%%%%%%%%%%%%%%%%%%%%%%%%%%%%

%
\title{\pier{Global solution and optimal control\\ of an epidemic propagation\\ with a heterogeneous diffusion}}
\author{}
\date{}
\maketitle
\Bcenter
\vskip-1.5cm
{\large\sc Pierluigi Colli$^{(1)}$}\\
{\normalsize e-mail: {\tt pierluigi.colli@unipv.it}}\\[.25cm]
{\large\sc Gianni Gilardi$^{(1)}$}\\
{\normalsize e-mail: {\tt gianni.gilardi@unipv.it}}\\[.25cm]
{\large\sc Gabriela Marinoschi$^{(2)}$}\\
{\normalsize e-mail: {\tt gabriela.marinoschi@acad.ro}}\\[.45cm]
$^{(1)}$
{\small Dipartimento di Matematica ``F. Casorati'', Universit\`a di Pavia}\\
{\small and Research Associate at the IMATI -- C.N.R. Pavia}\\ 
{\small via Ferrata 5, 27100 Pavia, Italy}\\[.2cm]
$^{(2)}$
{\small ``Gheorghe Mihoc-Caius Iacob'' Institute of Mathematical Statistics\\
and Applied Mathematics of the Romanian Academy}\\
{\small Calea 13 Septembrie 13, 050711 Bucharest, Romania}\\[.2cm]
\Ecenter

\Begin{abstract}
\noindent
\pier{In this paper, we explore the solvability and the optimal control problem for a
compartmental model based on reaction-diffusion partial differential
equations \gagi{describing a transmissible disease}. The nonlinear model takes into account the disease spreading due to the human social
diffusion, under a dynamic heterogeneity in infection risk. The analysis of the resulting system 
provides the existence proof for a global solution and \gagi{determines} the
conditions of optimality to reduce the concentration of the infected
population in certain spatial areas.}


\noindent {\bf AMS (MOS) Subject Classification:} 
35K55, % Nonlinear parabolic equations 
35K57, %	Reaction-diffusion equations  
35Q92, %	PDEs in connection with biology, chemistry and other natural sciences
46N60, % Applications of functional analysis in biology and other sciences
49J20, % Calculus of variations and optimal control; optimization;  Existence theories in calculus of variations and optimal control 
49J50, % Calculus of variations and optimal control; optimization;  Frechet and Gateaux differentiability in optimization	
49K20, % Calculus of variations and optimal control; optimization: Optimality conditions for problems involving partial differential equations 
92D30. %Epidemiology {For medical applications}
\End{abstract}
\salta
\pagestyle{myheadings}
\newcommand\testopari{\sc Colli \ --- \ Gilardi \ --- \ Marinoschi}
\newcommand\testodispari{{\sc \gagi{Control of epidemic propagation with heterogeneous diffusion}}}
\markboth{\testopari}{\testodispari}
\finqui
%
%%%%%%%%%%%%%%%%%%%%%%%%%%%%%%%%%
%% very beginning
%%%%%%%%%%%%%%%%%%%%%%%%%%%%%%%%%

\section{Introduction}
\label{Intro}
\setcounter{equation}{0}

\pier{The great impact on the development of human society that infectious
diseases can have requires prevention and control policies significant for
public health. The recent outbreak of infectious diseases, in particular
the COVID-19 pandemic, has highlighted an important role played by global
public surveillance systems and response which can lead to \gagi{reducing} their
effects on the socioeconomic activities and human health. In the past
decades, a number of mathematical models were developed to investigate
infectious disease evolution and the control of their spreading (an overview 
can be found in~\cite{ManOno, Song-Xiao}, see the references therein). In most of the cases, 
the models present in the literature are based on systems of ordinary
differential equations (ODEs) in time, describing the compartmental
evolution of various types of populations possibly evolving during an
epidemic.}

\pier{From the scientific contributions there is evidence showing that environmental heterogeneity
and human mobility have a significant impact on the spread of infectious
diseases (cf.,~e.g., \cite{Murray, Riley}); in addition, let us mention that
epidemic models accounting for spatial diffusion have been proposed and
investigated since long time ago (one may see \cite{Mottoni-Orlandi-Tesei, 
Webb, Fitzgibbon-1, Fitzgibbon-2}). An SIS functional
partial differential model cooperated with spatial heterogeneity and lag
effect of media impact was studied in~\cite{Song-Xiao}. A family of
epidemiological models, that extend the classic
Susceptible-Infectious-Recovered/Removed (SIR) model to account for dynamic
heterogeneity in infection risk, was proposed in~\cite{Berestycki}. The
family of models takes the form of a system of reaction-diffusion equations
for a given population structured by heterogeneous susceptibility to infection.
Recently, a new epidemic diffusion model with nonlinear transmission rates
and diffusion coefficients was introduced and tested in \cite{Viguerie-21, 
Viguerie-20}, while in \cite{Auricchio-23} the authors proved
well-posedness for an initial-boundary value problem associated to a variant
of the compartmental model for COVID-19 studied in \cite{Viguerie-21,
Viguerie-20}.}
\revis{Global existence and large time behavior of the solution to a system describing a 
spatio-temporal spread of an infectious disease are deduced in \cite[Section~3]{morgan} 
as an application of the results obtained in the same work~\cite{morgan}, 
for a more general model of semilinear reaction-diffusion-advection systems, 
whose coefficients satisfy general assumptions. In fact, in~\cite{morgan} the authors 
deal with unique, globally defined uniformly bounded weak solutions and provide a main 
issue showing that the quasi-positive systems that satisfy an intermediate sum condition 
automatically give rise to a new class of $L^{p}$-energy type functionals that allow to derive 
uniform a priori bounds.}

\pier{Since control measures are essential for a disease
mitigation, optimal control studies have been proposed in the literature,
considering generally the vaccine as a control variable (cf., e.g., \cite{parino,xu,vija}. However, other
precursory policies imposed to lower the contagion focus on the reduction of
the reproduction rate, in particular of the transmission rates and of the
population movement. Identification of such coefficients have been done in 
\cite{Giordano, GM-AMO, GM-DCDS}. Another approach of COVID-19 control has been introduced 
in the literature via the mean-field control model (see \cite{LLTLO, LLLO}),
where the aim is controlling the propagation of epidemics on a spatial domain. In the paper \cite{LLTLO}
 the control variable, the spatial velocity, was introduced for the classical disease models, such as the SIR 
model, and some numerical algorithms based on proximal
primal-dual methods were provided. The same method was used in \cite{LLLO} by choosing two controls 
for the pandemic: relocation of populations and distribution of vaccines.
\gagi{Optimization solutions in view of reducing the
number of infective individuals by using as controls the transmission rates
have been proposed in \cite{Medhaoui} and \cite{CGMR5}.}  Also, we quote the paper~\cite{dOnofrio} for 
the study of a more general
problem of optimally controlling an epidemic outbreak of a disease
structured by age since exposure, with the aid of two types of control
instruments, namely social distancing and vaccination. In \cite{dOnofrio} the aim is
minimizing the direct health cost of the epidemic arising from the overall
epidemic incidence, as well as the indirect epidemic cost, namely the
broader societal and economic cost due to social distancing and its impact
on the labor force and production, on overall social and relational
activities.}

\pier{In this paper, we explore an optimal control problem for a compartmental
model based on partial differential equations (PDEs), to account for the
disease spreading due to the human social diffusion under a dynamic
heterogeneity in infection risk and provide the conditions of optimality to
reduce the concentration of the infected population in certain spatial areas.}

\pier{The model we consider describes the evolution of an epidemic spreading in a
non-sedentary population}
\begin{align}
  & \dt s 
  + \betai s\, i
  + \betae s\, e 
  - \div (\kas \nabla s)
  - \gamma r
  = 0
   && \hbox{in $Q$}
  \label{Iprima}
  \\
  & \dt e
  - \betai s\, i
  - \betae s \, e 
  + \sigma e
  + \phie e
  - \div (\kae \nabla e)
  = 0
  &&\hbox{in $Q$}
  \label{Iseconda}
  \\
  & \dt i 
  + \phir i 
  - \div (\kai \nabla i)
  - \sigma e
  = 0
  &&\hbox{in $Q$}
  \label{Iterza}
  \\
  & \dt r 
  - \phir i
  - \phie e 
  - \div (\kar \nabla r)
  + \gamma r
  = 0
  && \hbox{in $Q$}
  \label{Iquarta}  
  \\
  & \dn s = \dn e = \dn i = \dn r = 0
  &&\hbox{on $\Sigma$}
  \label{Ibc}
  \\
  & \soluz(0) = (\sz,\ez,\iz,\rz)
  &&\hbox{in $\Omega$}
  \label{Icauchy}
\end{align}
\Accorpa\Ipbl Iprima Icauchy
\pier{where
\Beq
  Q := \Omega\times(0,T)
  \aand
  \Sigma := \Gamma \times (0,T)
  \label{defQS}
\Eeq
$\Omega$ being a bounded domain of $\erre^d$ (with \gagi{$d = 1,\, 2$} or 3) 
within which the population is located; here, $\Gamma$~denotes} the boundary 
of~$\Omega$ and $T>0$ is a fixed final time.
Moreover, in \eqref{Ibc}, $\dn$~denotes the derivative in the direction
 of the outward unit normal field $\nnn$ on the boundary.

\pier{Let us describe the meaning of the physical variables and
the coefficients that enter the system (1.1)--(1.6). The functions $s$, $e$, 
$i$ and $r$ represent the susceptible population, the exposed population,
the infected population, and the recovered population, respectively. 
Note that the the exposed population and the infected population
refer to the asymptomatic and symptomatic persons, respectively.
Hence, it turns out that, as observed in COVID-19 epidemic, the exposed may also
spread the disease.
Here we
neglect both the newborn and the deceased individuals after reaching the
maximum life, by assuming that these populations have small sizes and do not
contribute essentially to the epidemic transmission during the time period $%
(0,T)$. Thus, in our model the natality rate and the natural mortality rate
(which is not related to the disease) are zero. The nonnegative coefficients 
$\beta _{i}$ and $\beta _{e}$ depend on space, time and the total living
population%
\Beq
  n := s + e + i + r \,.
  \label{Idefn}  
\Eeq
We also assume that this epidemic evolves also by a diffusive process in
which healthy and infected individuals spread with different diffusion
coefficients $\kappa _{s},$ $\kappa _{e},$ $\kappa _{i}$ and $\kappa _{r}$
varying with space, time and the total living population.
The diffusion coefficients are assumed to be bounded from below and from above by positive constants,
so that the above system enjoys a parabolic character.
Finally, $\phir$, $\phie$ and $\sigma$ are fixed positive constants and $\gamma$ is a nonnegative function that depends just on time.}

\pier{For the system \Ipbl\ we can prove} \revis{a general existence result} \pier{by exploiting the 
analysis performed in \cite{CGMR5} and succesfully applying the Schauder fixed point theorem. 
Further, we discuss two uniqueness and continuous dependence \gagi{results} in a reduced setting. 
Moreover, in this paper the interest is particularly
focused on the control of the individual movements,
characterized by the time and space dependent controls given by the diffusion
coefficients $\kappa _{s},$ $\kappa _{e},$ $\kappa _{i}$ and $\kappa _{r}$, in order
to maintain the infected population (both~$e$ and~$i$) at a minimum level in a
certain subdomain $\Omega _{C}\subset \Omega .$ More exactly, having the map
of the initial values of the infected populations in the disjunct subdomains 
$\Omega _{j}$ of $\Omega $, we want to control the diffusion of the infected
individuals from these subdomains in the time interval $(0,T),$ such \gagi{that, by
limiting the individual movements by optimal policies,} to preserve a certain
domain with a lower density of infected. To this end, we consider that 
$ \Omega $ is represented as a reunion of \gagi{subdomains $\Omega_{j}$, $j=1,...,m$,
which can be referred as geographic areas, such that
\begin{align}
&\Omega =\bigcup\limits_{j=1}^{m}\Omega_{j},\quad\  \Omega_{j}\cap 
\Omega_{k}=\varnothing \quad \mbox{if  }\ j,k\in\{1,...,m\} , \, j\not=k, 
\non
\\
&\quad \Omega_{j}\,  \hbox{ is measurable and } \, |\Omega_j| > 0 \, \hbox{ for } \, j=1,...,m.
\label{pier11} 
\end{align}
Hence, we aim to} control the diffusion within the cylindrical domains $Q_{j}=\Omega _{j}\times (0,T),$
$j=1,...,m,$ via the controls $\kappa _{s},$ $\kappa _{e},$ $\kappa _{i},$ $%
\kappa _{r}$ that are assumed to have constant values in the domains $Q_{j},$
represented by 
\begin{align}
&\kappa _{s}(x,t) =\sum\limits_{j=1}^{m}u^s_j\,\chi _{Q_{j}}(x,t),\quad 
\kappa _{e}(x,t)=\sum\limits_{j=1}^{m}u^e_j\,\chi _{Q_{j}}(x,t), \non\\
&\quad \kappa _{i}(x,t) =\sum\limits_{j=1}^{m}u^i_j\,\chi _{Q_{j}}(x,t),\quad
\kappa _{r}(x,t)=\sum\limits_{j=1}^{m}u^r_j\, \chi _{Q_{j}}(x,t), \quad (x,t)\in Q \label{pier3}
\end{align}%
where $\chi _{Q_{j}}$ is the characteristic function of $Q_{j},$ $j=1,...,m.$
Obviously, the controls are the scalars $u^s_j,$ $u^e_j,$ $u^i_j,$
$u^r_j,$ but, for the sake of brevity, in most of our calculations we shall refer to the functions 
$\kappa_{s},$ $\kappa_{e},$ $\kappa _{i},$ $\kappa _{r}.$}

\pier{Then, assuming that $\Omega _{C}$ is a measurable set and denoting $Q_{C}=\Omega
_{C}\times (0,T),$ for a given positive coefficient $\alpha$ we introduce the cost functional 
\begin{align}
   \calJ(\kas,\kae,\kai,\kar,e,i)
  = \frac 12 \intQC \bigl( |e|^2 + |i|^2 \bigr)
  + \frac \alpha 2 \intQ \bigl( |\kas|^2 + |\kae|^2 + |\kai|^2 + |\kar|^2 \bigr)
  \label{Icost}
\end{align}
and state the control problem by 
\begin{align*}
   \hbox{minimize \ $\calJ(\kas,\kae,\kai,\kar,e,i)$ \ under the constraint }\ 
\control \in \Uad
 % \label{IcontrolPbl}
\end{align*}
subject to \Ipbl, where 
\Beq 
\Uad=\Uad^s\times\Uad^e\times\Uad^i\times\Uad^r
\label{pier2}
\Eeq
and each factor 
(which is the set of the admissible components of the controls) has the structure
\Beq
  \Uad^s := \biggl\{ \kas:=\somma j1{m} u_j^s \, \chi_{Q_j} :\ u_j^s\in[\umin^{s,j},\umax^{s,j}] \ \hbox{for $j=1,\dots,m$} \biggr\}
  \label{IdefUads}
\Eeq
for some fixed closed intervals $[\umin^{s,j},\umax^{s,j}]$, $j=1,\dots,m$. Here, only the first 
set $\Uad^s$ is precisely defined, but $\Uad^e, \, \Uad^i, \,\Uad^r$ are completely analogous.
\gagi{Next, for the results obtained for this optimal control problem, we refer the reader to the sequel and, in particular, to Section~\ref{CONTROL}.}} 

\pier{The paper is organized as follows. 
In the next section, we list our assumptions and notations
and state our results.
The proofs of our results regarding the existence and uniqueness for the state system \Ipbl\ and the continuous dependence of the solution on the controls are given in Section~\ref{WELLPOSEDNESS}, 
while Section~\ref{CONTROL} is devoted to the study of the optimal control problem.}

%%%%%%%%%%%%%%%%%%%%%%%%%%%%%%%%%%%%%%%%%%%%%%%%%%%%%%%%%%%%%%%%%%%%%%%%

\section{Statement of the problem and results}
\label{STATEMENT}
\setcounter{equation}{0}

In this section, we state precise assumptions and present our results.
First of all, the set $\Omega\subset\erre^d$, $d=1,2,3$, 
is~assumed to be bounded, connected and smooth.
Next, if $X$ is a Banach space, then $\norma\cpto_X$ denotes its norm,
with the only exception of the space $H$ defined below
and the $L^\infty$ spaces constructed on $\Omega$, $(0,T)$ and~$Q$, 
whose norms will be indicated by $\norma\cpto$ (i.e., without any subscript)
and by~$\norma\cpto_\infty\,$, respectively. 
Moreover, for simplicity, we use the same symbol for the norm in $X$ 
and that in any power of~$X$.
We also introduce 
\Beq
  H := \Ldue
  \aand 
  V := \Huno
  \label{defspazi}
\Eeq
\pier{and adopt} the framework of the Hilbert triplet $(V,H,\Vp)$
obtained by identifying $H$ with a subspace of the dual space $\Vp$ of $V$ in the usual way, namely, in order 
that $\<z,v>=\iO zv$ for every $z\in H$ and $v\in V$,
where $\<\cpto,\cpto>$ is the duality pairing between $\Vp$ and~$V$. 

\vskip 2mm

Now, we list the assumptions on the structure of the system.
We recall \eqref{defQS} for the definition of~$Q$
and we postulate~that
\begin{align}
  & \kas \,,\, \kae \,,\, \kai \,,\, \kar : Q \times \erre \to \erre
  \quad \hbox{are \Carath\ functions satisfying}
  \non
  \\
  & \quad \kamin \leq \kas(x,t,\zeta) \,,\, \kae(x,t,\zeta) \,,\, \kai(x,t,\zeta) \,,\, \kar(x,t,\zeta) \leq \kamax
  \non
  \\
  & \quad \hbox{\aaQ, every $\zeta\in\erre$ and some positive constants $\kamin$ and $\kamax$}
  \label{hpk}
  \\
  & \betai \,,\, \betae : Q \times \erre \to \erre
  \quad \hbox{are \Carath\ functions satisfying}
  \non
  \\
  & \quad 0 \leq \betai(x,t,\zeta) \,,\, \betae(x,t,\zeta) \leq \betamax
  \non
  \\
  & \quad \hbox{\aaQ, every $\zeta\in\erre$ and some positive constant $\betamax$}
  \label{hpbeta}
  \\
  & \hbox{$\phie$, $\phir$ and $\sigma$ are positive constants}
  \label{hpcoeff}
  \\
  & \gamma \in L^\infty(0,T) \quad \hbox{satisfies} \quad
  0 \leq \gamma(t) \leq \gammamax
  \non
  \\
  & \quad \hbox{\aat\ and some positive constant $\gammamax$}.
  \label{hpgamma}
\end{align}
\Accorpa\HPstruttura hpk hpgamma

\Brem
\label{Carath}
We recall that \gagi{$f:Q\times\erre\to\erre$ is a \Carath\ function if
\begin{align*}
&\hbox{$(x,t)\mapsto f(x,t,\zeta)$ is measurable on $Q$ for every $\zeta\in\erre$,}
\\[1mm]
&\hbox{$\zeta\mapsto f(x,t,\zeta)$ is continuous on $\erre$ \aaQ.}
\end{align*}
If} $f:Q\times\erre$ is a \Carath\ function
we still term $f$ the corresponding Nemytskii operator in the space of measurable functions on~$Q$, i.e.,
if $v:Q\to\erre$ is measurable, we use the following abbreviation
\Beq
  \hbox{$f(v)$ \ denotes the function \ $Q\mapsto\erre$ \ given by} \quad
  (x,t) \mapsto f(x,t,v(x,t)) \,.
  \label{abbrev}
\Eeq
We notice that
\Beq
  v_k \to v \quad \aeQ
  \quad \hbox{implies that} \quad
  f(v_k) \to f(v) \quad \aeQ .
  \label{carath}
\Eeq
\Erem

For the initial data, we postulate~that
\Beq
  \sz \,,\, \ez \,,\,\iz \,,\,\rz \in \Linfty
  \quad \hbox{are nonnegative}.
  \label{hpz}
\Eeq
Then, the state system related to the control problem we want to discuss is the following:
we look for a quadruplet $\soluz$ enjoying the regularity properties
\begin{align}
  & s ,\, e ,\, i ,\, r \in \H1\Vp \cap \L2V \emb \C0H
  \label{regsoluz}
  \\
  & s ,\, e ,\, i ,\, r \geq 0
  \quad \aeQ
  \label{possoluz}
  \\
  & s ,\, e ,\, i ,\, r \in \LQ\infty
  \label{bddsoluz}
\end{align}
\Accorpa\Regsoluz regsoluz bddsoluz
and satisfying the variational equations
\begin{align}
  & \< \dt s , v >
  + \iO \bigl( \betai(n) \, s i + \betae(n) \, s e \bigr) \, v
  + \iO \kas(n) \nabla s \cdot \nabla v
  - \iO \gamma r \, v
  = 0
  \label{prima}
  \\
  & \< \dt e , v >
  - \iO \bigl( \betai(n) \, s i + \betae(n) \, s e \bigr) \, v
  + \iO (\sigma+\phie) e \, v
  + \iO \kae(n) \nabla e \cdot \nabla v
  = 0
  \label{seconda}
  \\
  & \< \dt i , v >
  + \iO \phir \, i \, v
  + \iO \kai(n) \nabla i \cdot \nabla v
  - \iO \sigma e \, v
  = 0
  \label{terza}
  \\
  & \< \dt r , v >
  - \iO \bigl( \phir i + \phie e \bigr) \, v
  + \iO \kar(n) \nabla r \cdot \nabla v
  + \iO \gamma r \, v
  = 0
  \label{quarta}
\end{align}
\aet\ and for every $v\in V$, where \pier{(cf.~\eqref{Idefn})}
\Beq
  n := s + e + i + r 
  \label{defn}
\Eeq
as well as the initial condition
\Beq
  \soluz(0) = (\sz,\ez,\iz,\rz) \,.
  \label{cauchy}
\Eeq
\Accorpa\Pbl prima cauchy

Our first result regards the existence of a solution to this problem.

\Bthm
\label{Existence}
Assume \HPstruttura\ on the structure of the system and \eqref{hpz} on the initial data.
Then, there exists at least a quadruplet $\soluz$ satisfying \Regsoluz\
that solves problem \Pbl and satisfies the stability estimate
\Beq
  \norma\soluz_{\H1\Vp\cap\C0H\cap\L2V\cap\LQ\infty}
  \leq K_1
  \label{stab}
\Eeq
with some positive constant $K_1>0$ that depends only on
$\Omega$, $T$, the constants $\kamin$, $\kamax$, $\betamax$, $\gammamax$, $\phie$, $\phir$ and $\sigma$,
and the initial data.
\Ethm

We cannot prove uniqueness in the above general setting, unfortunately,
and we are able to show uniqueness only under suitable assumptions.
Namely, we first suppose~that
\Beq
  \kas = \kae = \kai = \kar
  \quad \hbox{do not explicitely depend on $(x,t)$}.
  \label{hpuniqueness}
\Eeq

\Bthm
\label{Uniqueness}
Besides the assumptions of Theorem~\ref{Existence}, assume \eqref{hpuniqueness}.
Then, problem \Pbl\ has a unique solution.
\Ethm

Under different assumptions we can prove both uniqueness and continuous dependence of the solution on the diffusion coefficients.
We stress that we do not need any longer that these coefficients are the same.
We recall that we have termed $\zeta$ the last variable in the assumptions regarding coefficients and suppose~that
\begin{align}
  & \hbox{the diffusion coefficients $\kas$, $\kae$, $\kai$ and $\kar$ do not depend on $\zeta$}
  \label{hpdiffusion}
  \\
  & |\betai(x,t,\zeta_1) - \betai(x,t,\zeta_2)|
  + |\betae(x,t,\zeta_1) - \betae(x,t,\zeta_2)|
  \leq L \, |\zeta_1-\zeta_2|
  \non
  \\
  & \quad \hbox{\aaQ, every $\zeta_1,\zeta_2\in\erre$ and some positive constant $L$}.
  \label{hpbetalip}
\end{align}
We will need a reinforcement of the above assumption only later~on 
(when dealing with the contol problem).
Nevertheless, we prefer to state it here.
We use the notation 
\Beq
  \hbox{$\betai'$ and $\betae'$ denote the derivatives with respect to $\zeta$ of $\betai$ and $\betae$, respectively}
  \non
\Eeq
(which exist almost everywhere due to \Lip\ continuity)
and assume that
\Beq
  \hbox{$\betai'$ and $\betae'$ are \Carath\ functions}.
  \label{hpbetaprimo}
\Eeq
For clarity, we write the form that problem \Pbl\ takes under assumption \eqref{hpdiffusion}:
the diffusion coefficients are just functions on $Q$ and $n$ enters just the coefficients $\betae$ and~$\betai$.
\begin{align}
  & \< \dt s , v >
  + \iO \bigl( \betai(n) \, s i + \betae(n) \, s e \bigr) \, v
  + \iO \kas \nabla s \cdot \nabla v
  - \iO \gamma r \, v
  = 0
  \label{primabis}
  \\
  & \< \dt e , v >
  - \iO \bigl( \betai(n) \, s i + \betae(n) \, s e \bigr) \, v
  + \iO (\sigma+\phie) e \, v
  + \iO \kae \nabla e \cdot \nabla v
  = 0
  \label{secondabis}
  \\
  & \< \dt i , v >
  + \iO \phir \, i \, v
  + \iO \kai \nabla i \cdot \nabla v
  - \iO \sigma e \, v
  = 0
  \label{terzabis}
  \\
  & \< \dt r , v >
  - \iO \bigl( \phir i + \phie e \bigr) \, v
  + \iO \kar \nabla r \cdot \nabla v
  + \iO \gamma r \, v
  = 0
  \label{quartabis}
  \\
  & \soluz(0) = (\sz,\ez,\iz,\rz) 
  \label{cauchybis}
\end{align}
where $n:=s+e+i+r$ as before.
\Accorpa\Pblbis primabis cauchybis

\Bthm
\label{Contdep}
Besides the assumptions of Theorem~\ref{Existence}, assume \accorpa{hpdiffusion}{hpbetalip}.
Then, problem \Pblbis\ has a unique solution.
Moreover, if $\kas_j$, $\kae_j$, $\kae_j$ and $\kar_j$, $j=1,2$,
are two choices of $\kas$, $\kae$, $\kai$ and~$\kar$, respectively,
and $(s_j,e_j,i_j,r_j)$ are the corresponding solutions,
then the estimate
\begin{align}
  & \norma{(s_1,e_1,i_1,r_1) - (s_2,e_2,i_2,r_2)}_{\H1\Vp\cap\C0H\cap\L2V}
  \non
  \\
  & \leq K_2 \, \norma{(\kas_1,\kae_1,\kai_1,\kar_1) - (\kas_2,\kae_2,\kai_2,\kar_2)}_\infty
  \label{contdep}
\end{align}
holds true with some positive constant $K_2>0$ that depends only on
$\Omega$, $T$, the constants $\kamin$, $\kamax$, $\betamax$, $\gammamax$, $\phie$, $\phir$ and $\sigma$,
the \Lip\ constant~$L$ and the initial data.
\Ethm

The above result is crucial for the the control problem we introduce at once.
The controls are the diffusion coefficients and we aim at controlling the exposed population and the infected population
in a region $\OmegaC\subset\Omega$.
Our theory needs the framework of Theorem~\ref{Contdep}.
Besides this, we assume~that
\begin{align}
  & \hbox{$\OmegaC$ \ is a measurable subset of $\Omega$ \ and \ $\QC:=\OmegaC\times(0,T)$}
  \label{hpOmegaC}
  \\
  & \hbox{$\alpha$ \ is a positive constant}.
  \label{hpalpha}
\end{align}
Then, the cost functional is \pier{defined by \eqref{Icost} for  $(\kas,\kae,\kai,\kar,e,i) \in (\LQ2)^6$, in principle. However, in view \gagi{of~\eqref{pier11}, \eqref{pier2} and~\eqref{IdefUads}},  the set $\Uad $ of admissible controls 
satisfies the general properties}
\begin{align}
%  & \calJ(\kas,\kae,\kai,\kar,e,i)
%  := \frac 12 \intQC \bigl( |e|^2 + |i|^2 \bigr)
%  + \frac \alpha 2 \intQ \bigl( |\kas|^2 + |\kae|^2 + |\kai|^2 + |\kar|^2 \bigr)
%  \label{cost}
%  \\
  & \hbox{$\Uad\subset(\LQ\infty)^4$ \ is convex and closed}
  \label{hpUad}
  \\
  & \kamin \leq \kas \,,\, \kae \,,\, \kai \,,\, \kar \leq \kamax
  \quad \hbox{\aeQ, \ for every $\control\in\Uad$}
  \non
  \\
  & \quad \hbox{and some positive constants $\kamin$ and $\kamax$}
  \label{hpkaminmax}
  \\
  & \hbox{the subspace $\Span(\Uad)$ is finite dimensional} 
  \label{hpfindim}
\end{align}
\Accorpa\HPcontrol hpOmegaC hpfindim
\pier{provided that
\begin{align}
&\hbox{the variability intervals 
$[\umin^{s,j},\umax^{s,j}]$, $[\umin^{e,j},\umax^{e,j}]$,}
\non 
\\
&\quad\hbox{ $[\umin^{i,j},\umax^{i,j}]$, $[\umin^{r,j},\umax^{r,j}]$, 
$j=1,\dots,m$, are  all contained in $[\kamin,\kamax]$.}
\label{pier10}
\end{align}
Then,} the control problem is the following
\begin{align}
  & \pier{\hbox{minimize \ $\calJ(\kas,\kae,\kai,\kar,e,i)$, \ for $\control \in \Uad\, $, \ where}}
%  \non
%  \\
%  & \quad \control \in \Uad
  \non
  \\
  & \quad \hbox{$e$ and $i$ are the \pier{components} of the solution to problem \Pblbis}
  \non
  \\
  & \quad \hbox{corresponding to the diffusion coefficients $\kas$, $\kae$, $\kai$ and $\kar$}.
  \label{controlPbl}
\end{align}

%\Brem
%\label{Findim}
%Assumption \eqref{hpfindim} is due to technical reasons (see the proof of the forthcoming Theorem~\ref{Optimum}).
%Nevertheless, we can describe a significant situation in which it is satisfied.
%Assume that $\Uad=\Uad^s\times\Uad^e\times\Uad^i\times\Uad^r$ and that each factor 
%(which is the set of the admissible components of the controls)
%has the structure
%that we describe just for the first one.
%Fix a measurable finite partition of unity~$\graffe{\chi_j:\ j=1,\dots,m_s}$, i.e.,
%each $\chi_j$ is a nonnegative function in $\LQ\infty$ and $\somma j1{m_s}\chi_j=1$ \aeQ.
%For, $j=1,\dots,m_s$, fix a closed interval $[\umin^{s,j},\umax^{s,j}]\subset[\kamin,\kamax]$.
%Then, set
%\Beq
%  \Uad^s := \Bigl\{ \kas:={\textstyle\somma j1{m_s} u_j^s \, \chi_j} :\ u_j^s\in[\umin^{s,j},\umax^{s,j}] \ \hbox{for $j=1,\dots,m_s$} \Bigr\}.
%  \label{defUads}
%\Eeq
%Then, it is clear that $\kamin\leq\kas\leq\kamax$ \aeQ\ for every $\kas\in\Uad^s$
%and that the ``true'' controls are given by the coefficients 
%that enter the definitions of $\Uad^s$, $\Uad^e$, $\Uad^i$ and~$\Uad^r$ (e.g., $u_j^s$ for $\Uad^s$).
%A particular but important case is the following:
%each $\chi_j$ is the characteristic function of some measurable subset $Q_j\subset Q$ of positive measure
%and $\graffe{Q_j}$ is a partition of~$Q$.
%For instance, $Q_j=\Omega_j\times(0,T)$ and $\graffe{\Omega_j}$ is a partition of~$\Omega$.
%\Erem

\pier{The existence of optimal controls and the necessary optimality conditions will be discussed in 
Section~\ref{CONTROL}.} We conclude this section by recalling the Young inequality 
\Beq
  a\,b \leq \delta\,a^2 + \frac 1{4\delta} \, b^2
  \quad \hbox{for all $a,b\in\erre$ and $\delta>0$}
  \label{young}
\Eeq
which will be repeatedly used throughout the paper
along with the \Holder\ and Schwarz inequalities,
and by stating a general rule concerning the constants. 
The small-case symbol $c$ stands for possibly different constants
(whose actual values may change from line to line and even within the same line)
that depend only on~$\Omega$, $T$, the structure of the system,
and the constants and the norms of the functions involved in the assumptions of the statements.
In particular, the values of $c$ are independent of the control variables $\kas$, $\kae$, $\kai$ and~$\kar$.
A~small-case symbol with a subscript like $\cdelta$
indicates that the constant may depend on the parameter~$\delta$, in addition.
On the contrary, we mark precise constants that we can refer~to
by using different symbols, e.g., capital letters.

%%%%%%%%%%%%%%%%%%%%%%%%%%%%%%%%%%%%%%%%%%%%%%%%%%%%%%%%%%%%%%%%%%%%%%%%

\section{Well-posedness}
\label{WELLPOSEDNESS}
\setcounter{equation}{0}

Well-posedenss in a particular case is proved in \cite[Thm~2.1]{CGMR5}.
We present this result in the form of a lemma.

\Blem
\label{ThmCGMR}
Assume \HPstruttura\ on the structure of the system 
and that all the coefficients $\kas$, $\kae$, $\kai$, $\kar$, $\betai$ and $\betae$
depend just on space and time and are independent of the last variable.
Moreover, assume that the initial data satisfy \eqref{hpz}.
Then, there exists a unique quadruplet $\soluz$ satisfying \Regsoluz\
and solving problem \Pbl.
Moreover, the stability estimate
\Beq
  \norma\soluz_{\H1\Vp\cap\C0H\cap\L2V\cap\LQ\infty}
  \leq K_0
  \label{stabCGMR}
\Eeq
holds true with some constant $K_0>0$ that depends only on
$\Omega$, $T$, the constants $\kamin$, $\kamax$, $\betamax$, $\gammamax$, $\phie$, $\phir$ and $\sigma$,
and the initial data.
\Elem

We use this result in our proofs of \pier{Theorems~\ref{Existence}} and~\ref{Uniqueness}.

\subsection{Existence}
\label{EXISTENCE}

We prove \pier{Theorem~\ref{Existence}} by applying the Schauder fixed point theorem.
Hence, we have to choose a proper subset $\calK$ of some Banach space $\calV$
that is convex and closed, and the continuous map $\Phi:\calK\to\calK$ 
in order that $\Phi(\calK)$ is relatively compact and a fixed point of $\Phi$ 
is a solution to the problem we want to solve.
To this end, we~set
\Beq
  \calV := \L2H
  \aand
  \calK := \graffe{v \in\calV:\ 0\leq v\leq 4\,K_0 \ \aeQ}
  \label{schauder}
\Eeq
where $K_0$ is the constant that appears in \eqref{stabCGMR}.
Then, we define $\Phi$ as follows.
For every $z\in\calK$, we consider the system given by the equations \accorpa{prima}{quarta}
where $n$ is replaced by $z$ in the argument of $\kas$, $\kae$, $\kai$, $\kar$, $\betai$ and~$\betae$,
complemented by the initial conditions~\eqref{cauchy}.
Then, we can apply the lemma: this system has a unique solution 
that we term $(s^z,e^z,i^z,r^z)$ and also satisfies~\eqref{stabCGMR}.
So, we~set
\Beq
  \Phi(z) := s^z + e^z + i^z + r^z \,.
  \label{defPhi}
\Eeq
Then, $\Phi(z)$ is nonnegative.
Moreover, $\norma{\Phi(z)}_\infty\leq\norma{s^z}_\infty+\norma{e^z}_\infty+\norma{i^z}_\infty+\norma{r^z}_\infty\leq 4\,K_0$.
Therefore $\Phi(z)\in\calK$.
Furthermore, it is clear that any fixed point of $\Phi$ is a solution to the original problem \Pbl.
To conclude the proof, it remains to show that $\Phi$ is continuous and that $\Phi(\calK)$ is relatively compact.
The latter immediately follows since \eqref{stabCGMR} applied to $(s^z,e^z,i^z,r^z)$ implies that
$\Phi(z)$ is bounded in $\H1\Vp\cap\L2V$,
so that we can apply the Aubin--Lions lemma (see, e.g., \cite[Thm.~5.1, p.~58]{Lions})
and conclude that $\Phi(\calK)$ is relatively compact.
Now, we prove that $\Phi$ is continuous.
To this end, we fix $z\in\calK$ and a sequence $\graffe{z_k}$ in $\calK$ that converges to $z$
and we show that $\Phi(z_k)$ converges to~$\Phi(z)$.
To simplify the notation, we write $(s_k,e_k,i_k,r_k)$ in place of $(s^{z_k},e^{z_k},i^{z_k},r^{z_k})$.
We also set $n_k:=\Phi(z_k)=s_k+e_k+i_k+r_k$.
By applying \eqref{stabCGMR} to $(s_k,e_k,i_k,r_k)$ and owing to standard compactness arguments, we deduce~that
\Beq
  (s_k,e_k,i_k,r_k) \to \soluz
  \label{convsoluzk}
\Eeq
weakly \pier{in $(\H1\Vp\cap\L2V)^4$ and weakly star in $\LQ\infty^4$},
for some limiting quadruplet $\soluz$ and for some subsequence.
We prove that $s+e+i+r$ coincides with $\Phi(z)$.
Once this is shown, since the limit has been characterized, we also conclude that the whole sequence $\Phi(z_k)$ converges to~$\Phi(z)$
and the proof is complete. 
To prove that $s+e+i+r$ coincides with $\Phi(z)$, is suffices to show that 
$\soluz$ is the solution to the system given by the equations \accorpa{prima}{quarta}
where $n$ is replaced by $z$ in the argument of $\kas$, $\kae$, $\kai$, $\kar$, $\betai$ and~$\betae$,
complemented by the initial conditions~\eqref{cauchy}.
The latter are clearly satisfied.
Let us discuss just the first equation since the others are analogous.
An equivalent formulation of our thesis is the following time-integrated version:
\begin{align}
  & \ioT \< \dt s(t) , v(t) > \, dt
  + \intQ \bigl( \betai(z) \, s i + \betae(z) \, s e \bigr) \, v
  + \intQ \kas(z) \nabla s \cdot \nabla v
  - \intQ \gamma r \, v
  = 0
  \non
  \\
  & \quad \hbox{for every $v\in \L2V$}.
  \label{peresistenza}
\end{align}
By the Aubin--Lions lemma, the convergence \eqref{convsoluzk} is also strong in $\L2H$ 
and, without loss of generality, almost every\pier{where} in~$Q$.
Then, by Remark~\ref{Carath}, we have that $\betai(z_k)$ converges to $\betai(z)$ \aeQ,
and analogous conclusions hold regarding $\betae$ and~$\kas$.
Since, \aeQ, it holds~that 
\Beq
  |\betai(z_k) s_k i_k| \leq \betamax K_0^2 \,, \quad 
  |\betae(z_k) s_k e_k| \leq \betamax K_0^2 
  \aand
  |\kas(z_k) \nabla v| \leq \kamax |\nabla v|
  \non
\Eeq
we can apply the Lebesgue dominated convergence theorem and deduce that
\begin{align}
  & \betai(z_k) s_k i_k \to \betai(z) s i 
  \aand
  \betae(z_k) s_k e_k \to \betae(z) s e
  \quad \hbox{strongly in $\L2H$}
  \non
  \\
  & \kas(z_k) \nabla v \to \kas(z) \nabla v
  \quad \hbox{strongly in $(\L2H)^d$}.
  \non
\end{align}
This, the weak convergence of $\dt s_k$, $\nabla s_k$ and $r_k$ given by \eqref{convsoluzk},
and the boundedness of~$\gamma$
ensure that we can pass to the limit in the analogue of \eqref{peresistenza} 
related to $z_k$ and $(s_k,e_k,i_k,r_k)$,~i.e.,
\Beq
  \ioT \< \dt s_k(t) , v(t) > \, dt
  + \intQ \bigl( \betai(z_k) \, s_k i_k + \betae(z_k) \, s_k e_k \bigr) \, v
  + \intQ \kas(z_k) \nabla s_k \cdot \nabla v
  - \intQ \gamma r_k \, v
  = 0
  \non
\Eeq
and obtain \eqref{peresistenza} itself.

%%%%%%%%%%%%%%%%%%%%%%%%%%%%%%%%%%%%%%%%%%%%%%%%%%%%%%%%%%%%%%%%%%%%%%%%

\subsection{Uniqueness and continuous dependence}
\label{CONTDEP}

This section is devoted to the proof of \pier{Theorems~\ref{Uniqueness}} and~\ref{Contdep}.

\step
Proof of \pier{Theorem~\ref{Uniqueness}}

We start with a preliminary observation.
The assumptions on the diffusion coefficients given in the statement
say that $\kas(n)$, $\kae(n)$, $\kai(n)$ and $\kar(n)$ are replaced by $\ka(n)$
where $\ka:\erre\to\erre$ is continuous and satisfies $\kamin\leq\ka(\zeta)\leq\kamax$ for every $\zeta\in\erre$.
Then, by adding all equations \accorpa{prima}{quarta} to each other, we obtain~that
\Beq
  \< \dt n , v >
  + \iO \nabla\Kappa(n) \cdot \nabla v
  = 0
  \quad \hbox{\aet, for every $v\in V$}
  \label{auxil}
\Eeq
where we have set
\Beq
  \Kappa(\zeta) := \int_0^\zeta \ka(\xi) \, d\xi
  \quad \hbox{for $\zeta\in\erre$}.
  \non
\Eeq
By taking $v=1$, we deduce that\pier{%
\Beq
\iO n(t) = \iO n(0) \quad \hbox{for every $t\in[0,T]$}.
\Eeq}%
From this, we derive a uniqueness result.
Assume that $n_1,n_2\in\H1\Vp\cap\L2V$ are two solutions to \eqref{auxil} with the same initial value.
Then $\iO(n_1-n_2)(t)=0$ for every $t\in[0,T]$.
By also recalling that $\Omega$ is connected, it \pier{is not difficult to show} that there exists a unique
$N\in\C0V$ that satisfies
\Beq
  \iO N(t) = 0
  \aand
  \iO \nabla N(t) \cdot \nabla v = \iO (n_1-n_2)(t) \, v
  \quad \hbox{for \pier{all} $t\in[0,T]$ and $v\in V$}.
  \non
\Eeq
Moreover, $N$ belongs to $\C0\Hdue$ since $n_1-n_2\in\C0H$ and $\Omega$ is smooth.
At this point, we write \eqref{auxil} for $n_1$ and~$n_2$,
take the difference and test it by~$N$.
Since $n_1-n_2=-\Delta N$, we obtain~that
\Beq
  \frac 12 \, \ddt \iO |\nabla N|^2
  + \iO \bigl( \Kappa(n_1) - \Kappa(n_2) \bigr) (n_1-n_2)
  = 0
  \quad \aet.
  \non
\Eeq
\pier{As} $\ka$ is positive, $\Kappa$ is \pier{an increasing function}, so that the above time derivative is nonpositive.
Since $n_i(0)=n_2(0)$ implies that $N(0)=0$, we conclude that $|\nabla N|^2$ vanishes identically.
Thus, $N$ is space independent.
\pier{Moreover, $N$ has also} zero mean value, \pier{and consequently} we conclude that $N=0$, whence $n_1=n_2$.

At this point, we can start our uniquenss proof.
We pick two solutions $(s_j,e_j,i_j,r_j)$, $j=1,2$. 
By applying the above observation to $n_j=s_j+e_j+i_j+r_j$,
we conclude that $n_1=n_2$.
It follows that $\kas(n_1)=\kas(n_2)$ and that similar equalities hold for 
$\kae$, $\kai$, $\kar$, $\betai$ and~$\betae$, i.e.,
these coeeficients can be considered as given functions on~$Q$.
Therefore, we can apply the uniqueness part of Lemma~\ref{ThmCGMR} 
and conclude that $(s_1,e_1,i_1,r_1)$ and $(s_2,e_2,i_2,r_2)$ coincide.

\step
Proof of Theorem~\ref{Contdep}

We pick two choices of the diffusion coefficients as in the statement
and a pair of arbitrary corresponding solutions.
We accordingly define $n_1$ and $n_2$ and set for brevity
\begin{align}
  & \kas := \kas_1 - \kas_2 \,, \quad
  \kae := \kae_1 - \kae_2 \,, \quad
  \kai := \kai_1 - \kai_2
  \aand
  \kar := \kar_1 - \kae_2
  \non
  \\
  & s := s_1 - s_2 \,, \quad
  e := e_1 - e_2 \,, \quad
  i := i_1 - i_2 \,, \quad
  r := r_1 - r_2 
  \aand
  n = n_1 - n_2 \,.
  \non
\end{align}
Then, we write equations \accorpa{primabis}{quartabis} for both coefficients and solutions,
take the differences, test them by $s$, $e$, $i$ and~$r$, respectively,
and rearrange a little.
We obtain~that
\begin{align}
  & \frac 12 \, \ddt \iO |s|^2                          % prima
  + \iO \kas_1 |\nabla s|^2
  = - \iO \kas \nabla s_2 \cdot \nabla s
  \non
  \\
  & \quad {}
  - \iO \bigl( \betai(n_1) - \betai(n_2) \bigr) s_1 \, i_1 \, s
  - \iO \betai(n_2) |s|^2 \, i_1
  - \iO \betai(n_2) s_2 \, i \, s
  \non
  \\
  & \quad {}
  - \iO \bigl( \betae(n_1) - \betae(n_2) \bigr) s_1 \, e_1 \, s
  - \iO \betae(n_2) |s|^2 \, e_1 
  - \iO \betae(n_2) s_2 \, e \, s
  + \iO \gamma r \, s \, \gagi,                        % fine prima
  \non
  \\
  \separa
  & \frac 12 \, \ddt \iO |e|^2                          % seconda
  + \iO \kae_1 |\nabla e|^2
  = - \iO \kae \nabla e_2 \cdot \nabla e
  \non
  \\
  & \quad {}
  + \iO \bigl( \betai(n_1) - \betai(n_2) \bigr) s_1 \, i_1 \, e
  + \iO \betai(n_2) s \, i_1 \, e
  + \iO \betai(n_2) s_2 \, i \, e
  \non
  \\
  & \quad {}
  + \iO \bigl( \betae(n_1) - \betae(n_2) \bigr) s_1 \, e_1 \, e
  + \iO \betae(n_2) s \, e_1 \, e
  + \iO \betae(n_2) s_2 \, |e|^2
  - (\sigma+\phie) \iO |e|^2 \, \gagi,                    % fine seconda
  \non
  \\
  \separa
  & \frac 12 \, \ddt \iO |i|^2                          % terza
  + \iO \kai_1 |\nabla i|^2
  = - \iO \kai \nabla i_2 \cdot \nabla i
  - \iO \phir \, |i|^2
  + \sigma \iO e \, i \, \gagi,                      % fine terza
  \non
  \\
  \separa
  & \frac 12 \, \ddt \iO |r|^2                          % quarta
  + \iO \kar_1 |\nabla r|^2
  = - \iO \kar \nabla r_2 \cdot \nabla r
  + \iO \bigl( \phir i + \phie e \bigr) \, r
  - \iO \gamma |r|^2 \,.                          % fine quarta
  \non
\end{align}
At this point, we add these equalities to each other.
The \lhs\ we obtain is bounded from below~by
\Beq
  \frac 12 \, \ddt \iO \bigl( |s|^2 + |e|^2 + |i|^2 + |r|^2 \bigr)
  + \kamin \iO \bigl( |\nabla s|^2 + |\nabla e|^2 + |\nabla i|^2 + |\nabla r|^2 \bigr)
  \non
\Eeq
and we have to estimate the \rhs\ from above.
For a while, the values of the constants denoted by $c$ can also depend on the solutions we have considered.
Some terms are nonpositive and can be neglected.
Let us treat just one of the terms involving gradients since the others are analogous.
By the Young inequality, we have \aet~that
\Beq
  - \iO \kas \nabla s_2 \cdot \nabla s
  \leq \norma\kas_\infty \iO |\nabla s_2| \, |\nabla s|
  \leq \frac{\pier{\kamin}} 2 \iO |\nabla s|^2
  + c \, \norma\kas_\infty^2 \, \normaV{s_2}^2 \,.
  \non
\Eeq
Since $s_2\in\L2V$, we deduce that
\Beq
  - \intQt \kas \nabla s_2 \cdot \nabla s
  \leq \frac{\pier{\kamin}} 2 \intQt |\nabla s|^2
  + c \, \norma\kas_\infty^2
  \non
\Eeq
where $Q_t:\Omega\times(0,t)$.
Let us now treat some of the terms involving products.
By \eqref{hpbetalip} and the \Holder\ and Young inequalities, we have \aet~that
\begin{align}
  & - \iO \bigl( \betai(n_1) - \betai(n_2) \bigr) s_1 \, i_1 \, s
  \leq L \, \norma{s_1}_\infty \, \norma{i_1}_\infty \iO |n| \, |s|
  \leq c \iO \bigl( |s|^2 + |e|^2 + |i|^2 + |r|^2 \bigr).
  \non
\end{align}
The next are even easier:
\begin{align}
  & - \iO \betai(n_2) |s|^2 \, i_1
  - \iO \betai(n_2) s_2 \, i \, s
  \leq \betamax \norma{i_1}_\infty \iO |s|^2
  + \betamax \norma{s_2}_\infty \iO |i| \, |s|
  \non
  \\
  & \leq c \iO \bigl( |s|^2 + |i|^2 \bigr) \,.
  \non
\end{align}
Since no new difficulty arises in estimating the other terms,
we quickly conclude.
After integrating with respect to time and owing to the initial conditions, 
we can apply the Gronwall lemma and obtain~that
\Beq
  \norma\soluz_{\C0H\cap\L2V}
  \leq c \, \norma{\control}_\infty \,.
  \label{quasicontdep}
\Eeq
Here, as said before, the value of $c$ can depend on the solutions we have considered.
Nevertheless, we can apply \eqref{quasicontdep} in the case $\kas=\kae=\kai=\kar=0$ and obtain that $\soluz=(0,0,0,0)$.
Since the solutions entering \eqref{quasicontdep} were arbitrary, this proves uniqueness.
Therefore, the solutions we have considered in the above proof are those given by our existence theorem 
and thus satisfy the stability estimate~\eqref{stab}.
Hence, the norms we have bounded by $c$ in the derivation 
can be bounded by the constant $K_1$ that appears in~\eqref{stab}.
Hence, \eqref{quasicontdep} actually holds with a constant $c$ 
that has the same dependence that the constant~$K_2$ of the statement has.
It remains to prove~that
\Beq
  \norma{\dt\soluz}_{\pier{\L2\Vp}}
  \leq c \, \norma{\control}_\infty 
  \label{tesidt}
\Eeq
with a similar constant~$c$.
We just consider the first component.
We come back to \eqref{primabis} written for both coefficients and solutions and take the difference.
Then, with an arbitrary $v\in\L2V$, we test the latter by~$v$ and have~that
\begin{align}
  & \< \dt s , v > 
  + \iO  \bigl(
    \betai(n_1) \, s_1 i_1
    - \betai(n_2) \, s_2 i_2
    + \betae(n_1) \, s_1 e_1
    - \betae(n_2) \, s_2 e_2
  \bigr) \, v
  \non
  \\
  & \quad {}
  + \iO (\kas_1 \nabla s_1 - \kas_2 \nabla s_2) \cdot \nabla v
  - \iO \gamma r \, v
  = 0 
  \quad \aet.
  \non 
\end{align}
By treating the products similarly as before, it is easy to prove~that
\Beq
  \< \dt s , v > 
  \leq c \iO (|s|+|e|+|i|+|r|) |v|
  + c \iO |\nabla s| \, |\nabla v|
  + c \, \norma \kas_\infty \iO |\nabla s_1| \, |\nabla v|.
  \non
\Eeq
We immediately deduce that
\Beq
  \< \dt s , v > 
  \leq c \, \bigl(
    \norma\soluz \, \norma v
    + \norma{\nabla s} \, \norma{\nabla v}
    + \norma\kas_\infty \norma{\nabla s_1} \, \norma{\nabla v} 
  \bigr) .
  \non
\Eeq
and by integrating over $(0,T)$ we infer that
\begin{align}
  & \ioT \< \dt s(t) , v(t) > \, dt
  \non
  \\
  & \leq c \, \bigl(
    \norma\soluz_{\L2H}
    + \norma s_{\L2V}
    + \norma\kas_\infty \norma{s_1}_{\L2V} 
  \bigr)  \norma v_{\L2V} \,.
  \non
\end{align}
By applying the stability estimate to $s_1$ and combining with \eqref{quasicontdep},
we conclude~that
\Beq
  \ioT \< \dt s(t) , v(t) > \, dt
  \leq c \, \norma\control_\infty \, \norma v_{\L2V} \,.
  \non
\Eeq
But this and the analogues for the other controls yield \eqref{tesidt} since $v$ is arbitrary in $\L2V$, 
and the proof is complete.

%%%%%%%%%%%%%%%%%%%%%%%%%%%%%%%%%%%%%%%%%%%%%%%%%%%%%%%%%%%%%%%%%%%%%%%%

\section{The control problem}
\label{CONTROL}
\setcounter{equation}{0}

This section is devoted to the study of the control problem \eqref{controlPbl} presented in the Introduction.
\pier{We first prove the existence of an optimal control.
Then, we derive a first order necessary condition for optimality.
Concerning the proof of this condition, we do not follow the standard technique
based on the \Frechet\ differentiability of the control-to-state mapping,
since the problem is particularly complicated.
So, we prefer to deduce the necessary condition
by directly acting on the first variation of the cost functional and proceeding in the direction of G\^ateax derivatives.
The condition we prove sounds as follows.
If $\controlstar$ is an optimal control and $\soluzstar$ is the corresponding optimal state, 
then there holds the variational inequality
\begin{align}
  & \intQ \bigl\{ 
    (\alpha \kasstar - \nabla\sstar \cdot \nabla p) (\kas-\kasstar) 
    + (\alpha \kaestar - \nabla\estar \cdot \nabla q) (\kae-\kaestar)
  \non 
  \\
  & \quad {}
    + (\alpha \kaistar - \nabla\istar \cdot \nabla w) (\kai-\kaistar) 
    + (\alpha \karstar - \nabla\rstar \cdot \nabla z) (\kar-\karstar) 
  \bigr\}
  \non
  \\
  & \geq 0
  \quad \hbox{for every $\control\in\Uad$}
  \non
\end{align}
where $\soluza$ is the solution to a proper adjoint problem (cf.~\eqref{primaa}--\eqref{cauchya}). 
%The latter is a backward parabolic problem that we do not describe here for brevity
%since the terms entering it are rather complicate.
%All the details are given in Section~\ref{CONTROL},
%where, in the forthcoming Remark~\ref{RemNC}, 
%we also specify the form that the above condition takes in the last situation of Remark~\ref{Findim},
%where the functions $\chi_j$ are orthogonal to each other with respect to the standard inner product of~$\LQ2$:
We will see that the optimal coefficients are just one-dimensional projections
of weighted mean values of suitable ingredients of the adjoint system
on the closed intervals that enter the particular definition of~$\Uad$.}

Besides the general assumptions \HPstruttura\ on the structure of the system and \eqref{hpz} on the initial data,
we need to use \accorpa{hpdiffusion}{hpbetalip} as well as \eqref{hpbetaprimo}.
We recall the form \Pblbis\ that the original state system takes under the above \pier{assumptions.
%and that the cost functional and the set of admissible controls 
%are ruled by \HPcontrol\ \pier{and \eqref{pier2}--\eqref{IdefUads}}.
It is understood that these assumptions are in force in the whole section
and we do not recall them in our statements.}

%%%%%%%%%%%%%%%%%%%%%%%%%%%%%%%%%%%%%%%%%%%%%%%%%%%%%%%%%%%%%%%%%%%%%%%%

\subsection{Existence of an optimal control}
\label{OPTIMUM}

\Bthm
\label{Optimum}
\pier{Assume \HPcontrol. Then} there exists an optimal control, i.e.,
there exist $\controlstar\in\Uad$ such that,
if $\control\in\Uad$ and $\soluzstar$ and $\soluz$ 
are the states corresponding to $\controlstar$ and $\control$, respectively,
then there holds the inequality
\Beq
  \calJ(\kasstar,\kaestar,\kaistar,\karstar,\estar,\istar)
  \leq \calJ(\kas,\kae,\kai,\kar,e,i) \,.
  \label{optimum}
\Eeq
\Ethm

\Bdim
We use the direct method.
Let $\Lam\geq0$ be the infimum of the cost functional subject to the constraints specified in \eqref{controlPbl},
let $\graffe{(\kas_k,\kae_k,\kai_k,\kar_k)}$ be a minimizing sequence 
and let $\graffe{(s_k,e_k,i_k,r_k)}$ be the sequence of the corresponding states.
Since $\Uad$ is bounded in $(\LQ\infty)^4$ and $\Span(\Uad)$ is finite dimensional, 
the sequence of controls has a strongly converging subsequence,
which we still term $\graffe{(\kas_k,\kae_k,\kai_k,\kar_k)}$ to simplify the notation.
By applying the stability estimate \eqref{stab} to the above sequence of states and owing to standard compactness results
(in~\gagi{particular} to the Aubin--Lions lemma),
we find a convergence subsequence, still termed $\graffe{(s_k,e_k,i_k,r_k)}$ to avoid a boring notation, in the related topology.
We thus have~that
\begin{align}
  & \calJ(\kas_k,\kae_k,\kai_k,\kar_k,e_k,i_k) \to \Lambda
  \label{convcosts}
  \\
  & (\kas_k,\kae_k,\kai_k,\kar_k) \to \controlstar
  \quad \hbox{strongly in $(\LQ\infty)^4$}
  \label{convcontrols}
  \\
  & (s_k,e_k,i_k,r_k) \to \soluzstar
  \non
  \\
  & \quad \hbox{weakly star in $(\H1\Vp\cap\L2V\cap\LQ\infty)^4$,}
  \non
  \\
  & \qquad \hbox{strongly in $(\L2H)^4$ and \aeQ}
  \label{convstates}
\end{align}
for some quadruplets $\controlstar$ and $\soluzstar$.
The former belongs to $\Uad$ since $\Uad$ is closed.
We now prove that $\soluzstar$ is the state corresponding to $\controlstar$, i.e.,
it solves \Pblbis\ where one reads $\controlstar$ in place of $\control$.
The initial condition \eqref{cauchybis} is satisfied since the sequence $\graffe{(s_k,e_k,i_k,r_k)}$ also converges weakly in $(\C0H)^4$.
We just prove that the first equation \eqref{primabis} is satisfied and we do not repeat this argument for the others.
Namely, we prove that the limiting quadruplet satisfies the equivalent time-integrated version of \eqref{primabis},~i.e.,
\begin{align}
  & \ioT \< \dt \sstar(t) , v(t) > \, dt
  + \intQ \bigl( \betai(\nstar) \, \sstar \istar + \betae(\nstar) \, \sstar \estar \bigr) \, v
  + \intQ \ka_s^*(\nstar) \nabla\sstar \cdot \nabla v
  \non
  \\
  & \quad {}
  - \intQ \gamma \rstar \, v
 = 0
  \quad \hbox{for every $v\in \L2V$}
  \label{peroptimum}
\end{align} 
where $\nstar:=\sstar+\estar+\istar+\rstar$.
Thus, we fix $v\in\L2V$.
By obviously definining $n_k:=s_k+e_k+i_k+r_k$, we have that
\Beq
  |\betai(n_k) \, s_k i_k v| \leq c \, |v|
  \aand
  |\betae(n_k) \, s_k e_k v| \leq c \, |v|
  \quad \aeQ
  \non
\Eeq
so that
\Beq
  \betai(n_k) \, s_k i_k v \to \betai(\nstar) \, \sstar \istar v
  \aand
  \betae(n_k) \, s_k e_k v\to \betae(\nstar) \, \sstar \estar v
  \non
\Eeq
strongly in $\L2H$ by the Lebesgue dominated convergence theorem.
On the other hand, \eqref{convcontrols} implies
\Beq
  \kas_k \nabla v \to \kasstar \nabla v
  \quad \hbox{strongly in $(\L2H)^d$}
  \non
\Eeq
and the analogues for the other controls.
Therefore, one can pass to the limit as $k$ tends to infinity in the analogue of \eqref{peroptimum}
written for $(\kas_k,\kae_k,\kai_k,\kar_k)$ and the corresponding states
and obtain \eqref{peroptimum} itself.
By strong convergence, we also have~that
\Beq
  \calJ(\kas_k,\kae_k,\kai_k,\kar_k,e_k,i_k) \to \calJ(\kasstar,\kaestar,\kaistar,\karstar,\estar,\istar)
  \non
\Eeq
whence $\calJ(\kasstar,\kaestar,\kaistar,\karstar,\estar,\istar)=\Lambda$.
We conclude that $\controlstar$ is an optimal control and $\soluzstar$ is the corresponding state.
\Edim

%%%%%%%%%%%%%%%%%%%%%%%%%%%%%%%%%%%%%%%%%%%%%%%%%%%%%%%%%%%%%%%%%%%%%%%%

\subsection{Necessary conditions for optimality}
\label{NC}

We fix an optimal control $\controlstar$ and the corresponding state $\soluzstar$.
Given any $\control\in\Uad$, we associate the following linearized system
\begin{align}
  & \< \dt\xi , v >
  + \iO \kasstar \nabla\xi \cdot \nabla v
  + \iO \bigl( A_1 \xi + B_1 \eta + C_1 \iota + D_1 \rho) v
  \non
  \\
  & = - \iO (\kas-\kasstar) \nabla\sstar \cdot \nabla v  \, \gagi,
  \label{primal}
  \\
  & \< \dt\eta , v >
  + \iO \kaestar \nabla\eta \cdot \nabla v
  + \iO \bigl( A_2 \xi + B_2 \eta + C_2 \iota + D_2 \rho) v
  \non
  \\
  & = - \iO (\kae-\kaestar) \nabla\estar \cdot \nabla v   \, \gagi,
  \label{secondal}
  \\
  & \< \dt\iota , v >
  + \iO \kaistar \nabla\iota \cdot \nabla v
  + \iO \bigl( A_3 \xi + B_3 \eta + C_3 \iota + D_3 \rho) v
  \non
  \\
  & = - \iO (\kai-\kaistar) \nabla\istar \cdot \nabla v    \, \gagi,
  \label{terzal}
  \\
  & \< \dt\rho , v >
  + \iO \karstar \nabla\rho \cdot \nabla v
  + \iO \bigl( A_4 \xi + B_4 \eta + C_4 \iota + D_4 \rho) v
  \non
  \\
  & = - \iO (\kar-\karstar) \nabla\rstar \cdot \nabla v    
  \label{quartal}
  \\
  & \soluzl(0) = (0,0,0,0) 
  \label{cauchyl}
\end{align}
\Accorpa\Pbll primal cauchyl
where we have~set
\begin{align}
  & A_1 := \betai'(\nstar) \sstar \istar
  + \betai(\nstar) \istar
  + \betae'(\nstar) \sstar \estar
  + \betae(\nstar) \estar
  \label{defA1}
  \\
  & B_1 := \betai'(\nstar) \sstar \istar
  + \betae(\nstar) \sstar
  + \betae'(\nstar) \sstar \estar
  \label{defB1}
  \\
  & C_1 := \betai'(\nstar) \sstar \istar
  + \betai(\nstar) \sstar
  + \betae'(\nstar) \sstar \estar
  \label{defC1}
  \\
  & D_1 := \betai'(\nstar) \sstar \istar
  + \betae'(\nstar) \sstar \estar
  - \gamma
  \label{defD1}
  \\
  & A_2 := - A_1 \,, \quad 
  B_2 := - A_2 + \sigma + \phie \,, \quad
  C_2 := - C_1 \,, \quad
  D_2 := - D_1 - \gamma
  \label{def2}
  \\
  & A_3 := 0 \,, \quad 
  B_3 := - \sigma \,, \quad
  C_3 := \phir \,, \quad
  D_3 := 0
  \label{def3}
  \\
  & A_4 := 0 \,, \quad 
  B_4 := - \phie \,, \quad
  C_4 := - \phir \,, \quad
  D_4 := \gamma \,.
  \label{def4}
\end{align}
We notice at once that problem \Pbll\ has a unique solution
\Beq
  \soluzl \in (\H1\Vp\cap\L2V)^4 
  \label{regsoluzl}
\Eeq
since it is a uniformly parabolic problem with bounded coefficients and known \rhs s belonging to $\L2\Vp$
if we read \accorpa{primal}{quartal} as an abstract equation in the framwork $(V,H,\Vp)^4$.
We have the following result:

\Bprop
\label{BadNC}
\pier{Assume \HPcontrol\ and} let $\controlstar$ and $\soluzstar$ be an optimal control and the corresponding state, respectively.
Then, for every $\control\in\Uad$ there holds the inequality
\begin{align}
  & \intQC (\eta \, \estar + \iota \, \istar)
  \non
  \\
  & \quad {}
  + \alpha \intQ \bigl( 
    (\kas-\kasstar) \, \kasstar
    + (\kae-\kaestar) \, \kasstar
    + (\kai-\kaistar) \, \kasstar
    + (\kar-\karstar) \, \kasstar
  \bigr) \geq 0
  \label{badNC}
\end{align}
where $\soluzl$ is the solution to the linearized system associated with $\control$.
\Eprop

\Bdim
We fix the admissible control $\control$ and the corresponding state $\soluz$ and,
for every $\lam\in(0,1)$, we introduce the incremented control
\Beq
  \controllam
  := \controlstar + \lam \bigl( \control - \controlstar \bigl)
  \label{controllam}
\Eeq
which belongs to $\Uad$ since $\Uad$ is convex.
We also~set
\Beq
  \hbox{$\soluzlam$ is the state corresponding to $\controllam$}
  \label{statelam}
\Eeq
and define the quotients
\Beq
  \xilam := \frac {\slam-\sstar} \lam \,, \quad 
  \etalam := \frac {\elam-\estar} \lam \,, \quad
  \iotalam := \frac {\ilam-\istar} \lam
  \aand
  \rholam := \frac {\rlam-\rstar} \lam \,.
  \label{quotients}
\Eeq
We prove that $(\xilam,\etalam,\iotalam,\rholam)$ converges in a suitable topology 
to the solution $\soluzl$ to the linearized system as $\lam$ tends to zero.
First, we notice the regularity
\Beq
  (\xilam,\etalam,\iotalam,\rholam) \in (\H1\Vp\cap\L2V)^4 \,.
  \label{regquotients}
\Eeq
Then, we write the system these quotients solve.
This is obtained as follows:
we write \Pblbis\ for the incremented control \eqref{controllam} and the corresponding state;
then we do the same for the optimal control and state;
finally, we take the differences and divide by~$\lambda$.
With the notations
\begin{align} 
  \betai^* := \betai(\nstar) 
  \aand
  \betai^\lam := \betai(\nlam)
  \quad \hbox{where} \quad
  \nlam := \slam + \elam + \ilam + \rlam
  \label{notationbeta}
\end{align}
and the analogous ones for $\betae^*$ and $\betae^\lam$,
and observing~that
\Beq
  \frac {\kaslam-\kasstar} \lam = \kas - \kasstar 
  \non
\Eeq
and that analogous identities hold for the other similar fractions,
we have~that
\begin{align}
  & \< \dt\xilam , v >
  + \iO \frac 1\lam (\betai^\lam-\betai^*) \, \slam \ilam v
  + \iO \betai^* \, (\xilam\ilam + \sstar\iotalam) \pier{v}
  \non
  \\
  & \quad {}
  + \iO \frac 1\lam (\betae^\lam-\betae^*) \, \slam \elam v
  + \iO \betae^* \, ( \xilam\elam + \sstar\etalam) v
  \non
  \\
  & \quad {}
  + \iO \kasstar \nabla\xilam \cdot \nabla v
  - \iO \gamma \rholam v
  = - \iO (\kas-\kasstar) \nabla\slam \cdot \nabla v \, \gagi,
  \label{primalam}
  \\
  \separa
  & \< \dt\etalam , v >
  - \iO \frac 1\lam (\betai^\lam-\betai^*) \, \slam \ilam v
  - \iO \betai^* \, (\xilam\ilam + \sstar\iotalam)\pier{v}
  \non
  \\
  & \quad {}
  - \iO \frac 1\lam (\betae^\lam-\betae^*) \, \slam \elam v
  - \iO \betae^* \, (\xilam\elam + \sstar\etalam) v
  \non
  \\
  & \quad {}
  + (\sigma+\phie) \iO \etalam v
  + \iO \kaestar \nabla\etalam \cdot \nabla v
  = - \iO (\kae-\kaestar)  \nabla\elam \cdot \nabla v \, \gagi,
  \label{secondalam}
  \\
  \separa
  & \< \dt\iotalam , v >
  + \phir \iO \iotalam v
  + \iO \kaistar \nabla\iotalam \cdot \nabla v
  \non
  \\
  & \quad {}
  - \sigma \iO \etalam v
  = - \iO (\kai-\kaistar)  \nabla\ilam \cdot \nabla v \, \gagi,
  \label{terzalam}
  \\
  \separa
  & \< \dt\rholam , v >
  - \phir \iO \iotalam v
  - \phie \iO \etalam v
  + \iO \gamma \rholam v
  + \iO \kaistar \nabla\iotalam \cdot \nabla v
  \non
  \\
  & = - \iO (\kar-\karstar)  \nabla\rlam \cdot \nabla v \, \gagi,
  \label{quartalam}
\end{align}
all the equations holding \aet\ and for every $v\in V$.
Moreover, the initial condition
\Beq
  (\xilam,\etalam,\iotalam,\rholam) = (0,0,0,0) 
  \label{cauchylam}
\Eeq
\Accorpa\Pbllam primalam cauchylam
is also satisfied.
In order to control the behavior as $\lambda$ tends to zero,
we perform an estimate.
We notice that the stability estimate \eqref{stab}
holds for both $\soluzstar$ and $\soluzlam$.
By applying the continuous dependence estimate~\eqref{contdep},
we have~that
\begin{align}
  & \norma{\soluzlam-\soluzstar}_{\H1\Vp\cap\L2V}
  \non
  \\
  & \leq c \, \norma{\controllam-\controlstar}_\infty
  \non
  \\
  & {} = c \, \lam \, \norma{\control-\controlstar}_\infty
  \leq c \, \lam
  \label{perstimalam}
\end{align}
whence
\Beq
  \norma{(\xilam,\etalam,\iotalam,\rholam)}_{\H1\Vp\cap\L2V} \leq c \,.
  \label{stimalam}
\Eeq
Therefore, thanks to standard compactness results, we have~that
\Beq
  (\xilam,\etalam,\iotalam,\rholam) \to \soluzl
  \quad \hbox{weakly in $(\H1\Vp\cap\L2V)^4$}
  \non
\Eeq
for some quadruplet $\soluzl$ satisfying \eqref{regsoluzl} 
(here and later on, for a subsequence, in principle;
however, after the proof we now perform,
the whole family converges, by uniqueness).
We prove that the limiting quadruplet solves problem \Pbll.
As usual, in letting $\lambda$ tend to zero,
we consider the time-integrated versions of both \Pbllam\ and \Pbll\ with time dependent test functions $v\in\L2V$.
For brevity, we do not treat all the terms of \Pbllam\ 
since many of them can be dealt with very easily, in particular since
\begin{align}
  & \soluzlam \to \soluzstar
  \aand
  (\xilam,\etalam,\iotalam,\rholam) \to \soluzl
  \non
  \\
  & \quad \hbox{strongly in $(\L2H)^4$ and \aeQ}
  \label{strongconvlam}
\end{align}
the former by \eqref{perstimalam} and the latter
thanks to the Aubin--Lions lemma (see, e.g., \cite[Thm.~5.1, p.~58]{Lions}).
Namely, we confine ourselves to consider the integral involving~$\betai^\lam$
(the same argument is used for~$\betae^\lam$), 
which is the most delicate.
For a fixed $v\in\L2V$, we prove~that
\Beq
  \intQ \frac 1\lam (\betai^\lam-\betai^*) \, \slam \ilam v
  \quad\to\quad \intQ \betai'(\nstar) (\xi+\eta+\iota+\rho) \, \sstar \istar v 
  \label{limprimotermine}
\Eeq
as $\lambda$ tends to zero.
We have \aeQ~that
\begin{align}
  & \frac {\betai^\lam-\betai^*} \lam
  = \frac 1\lam \int_{\nstar}^{\nlam} \betai'(\zeta) \, d\zeta
  = \frac 1\lam \int_0^1 \betai'(\nstar+\tau(\nlam-\nstar)) (\nlam-\nstar) \, d\tau
  \non
  \\
  & = \int_0^1 \betai'(\nstar+\tau(\nlam-\nstar)) \, d\tau \, (\xilam+\etalam+\iotalam+\rholam).
  \non
\end{align}
On the other hand, \eqref{strongconvlam} implies that $\nlam$ converges to $\nstar$ \aeQ.
Since $\betai'$ is a \Carath\ function (cf.~\eqref{hpbetaprimo}) we have that 
$\betai'(\nstar+\tau(\nlam-\nstar))$ converges to $\betai'(\nstar)$ \aeQ\ for every $\tau\in(0,1)$.
Finally, $\betai'$ is bounded by the \Lip\ constant~$L$ (see \eqref{hpbetalip}).
Hence, we have~that
\Beq
  \int_0^1 \betai'(\nstar+\tau(\nlam-\nstar)) \, d\tau
  \quad\to\quad \int_0^1 \betai'(\nstar) \, d\tau
  = \betai'(\nstar)
  \non
\Eeq
by the Lebesgue dominated convergence theorem.
Therefore
\Beq
  \frac 1\lam (\betai^\lam-\betai^*) \slam \ilam
  \quad\to\quad \betai'(\nstar) (\xi+\eta+\iota+\rho) \sstar \istar
  \quad \aeQ
  \non
\Eeq
as $\lambda$ tends to zero.
On the other hand, on account of \eqref{stimalam}, we also have~that
\Beq
  \norma{(1/\lam)(\betai^\lam-\betai^*)}_{\L2H}
  \leq (L/\lam) \, \norma{\nlam-\nstar}_{\L2H}
  = L \norma{\xilam+\etalam+\iotalam+\rholam}_{\L2H}
  \leq c 
  \non
\Eeq
whence also
\Beq
  \norma{(1/\lam)(\betai^\lam-\betai^*)\slam\ilam}_{\L2H}
  \leq c \,.
  \non
\Eeq
Hence
\Beq
  \frac 1\lam (\betai^\lam-\betai^*) \slam \ilam
  \quad\to\quad \betai'(\nstar) (\xi+\eta+\iota+\rho) \sstar \istar
  \quad \hbox{weakly in $\L2H$}
  \non
\Eeq
and \eqref{limprimotermine} follows.
This concludes the proof that $\soluzl$ is the solution to \Pbll.

At this point, we are ready to prove \eqref{badNC}.
Due to optimality, we have~that
\Beq
  \frac{\calJ(\kaslam,\kaelam,\kailam,\karlam,\elam,\ilam)
    - \calJ(\kasstar,\kaestar,\kaistar,\karstar,\estar,\istar)} \lam
  \geq 0
  \quad \hbox{for every $\lam\in(0,1)$}
  \label{perNC}
\Eeq
and we aim at letting $\lam$ tend to zero in this inequality.
We just consider two of the terms involved in~\eqref{perNC}, namely
\Beq
  \frac 12 \intQC \frac {|\elam|^2-|\estar|^2} \lam
  \aand
  \frac \alpha 2 \intQ \frac {|\kaslam|^2-|\kasstar|^2} \lam
  \non
\Eeq
since the others are analogous.
We have that
\Beq
  \frac 12 \intQC \frac {|\elam|^2-|\estar|^2} \lam
  = \intQC \frac {\elam-\estar} \lam \, \frac {\elam+\estar} 2
  = \intQC \etalam \, \frac {\elam+\estar} 2
  \non
\Eeq
and this converges to
\Beq
  \intQC \eta \, \estar
  \non
\Eeq
as $\lam$ tends to zero.
Similarly
\Beq
  \frac \alpha 2 \intQ \frac {|\kaslam|^2-|\kasstar|^2} \lam
  = \alpha \intQ \frac {\kaslam-\kasstar} \lam \, \frac {\kaslam+\kasstar} 2
  = \alpha \intQ (\kas-\kasstar) \, \frac {\kaslam+\kasstar} 2
  \non
\Eeq
converges to
\Beq
  \alpha \intQ (\kas-\kasstar) \, \kasstar \,.
  \non
\Eeq
Hence, \eqref{badNC} immediately follows.
\Edim

The result just proved is not satisfactory.
Indeed, the linearized problem \Pbll\ is involved infinitely many times
since $\control$ is arbitrary in~$\Uad$.
This trouble is bypassed by the introduction of a proper adjoint problem.
This is formally obtained by testing the equations \accorpa{primal}{quartal} by $p$, $q$, $z$ and~$w$, respectively,
integrating by parts in time and collecting the terms involving $\xi$, $\eta$, $\iota$ and~$\rho$, respectively.
It is the following
\begin{align}
  & - \< \dt p , v >
  + \iO \kasstar \nabla p \cdot \nabla v
  + \iO ( A_1 p + A_2 q + A_3 w + A_4 z ) v
  = 0
  \label{primaa}
  \\
  & - \< \dt q , v >
  + \iO \kaestar \nabla q \cdot \nabla v
  + \iO ( B_1 p + B_2 q + B_3 w + B_4 z ) v
  = \iOC \estar v
  \label{secondaa}
  \\
  & - \< \dt w , v >
  + \iO \kaistar \nabla w \cdot \nabla v
  + \iO ( C_1 p + C_2 q + C_3 w + C_4 z ) v
  = \iOC \istar v
  \label{terzaa}
  \\
  & - \< \dt z , v >
  + \iO \karstar \nabla z \cdot \nabla v
  + \iO ( D_1 p + D_2 q + D_3 w + D_4 z ) v
  = 0
  \label{quartaa}
  \\
  & (p,q,w,z)(T) = (0,0,0,0) 
  \label{cauchya}
\end{align}
\Accorpa\Pbla primaa cauchya
all the equations holding \aet\ and for every $v\in V$,
where $A_j$, $B_j$, $C_j$ and $D_j$, $j=1,\dots,4$, are defined in \accorpa{defA1}{def4}.
This is a backward linear parabolic problem with bounded coefficients.
Thus, it has a unique solution
\Beq
  (p,q,rw,z) \in (\H1\Vp\cap\L2V)^4
  \label{regsoluza}
\Eeq
which actually is more regular.
A~satisfactory \pier{and explicit} necessary condition for optimality is given in the following theorem,
which \pier{is our} last result. The precise form of the set $\Uad$ specified in the Introduction
is assumed. 

\Bthm
\label{GoodNC}
\pier{Assume \HPcontrol\ and \gagi{\eqref{pier11},} \eqref{pier2}, \eqref{IdefUads}, \eqref{pier10}. Let
$\controlstar$ be an optimal control of the form 
\begin{align}
&\kappa_{s}(x,t) =\sum\limits_{j=1}^{m}u^*_{s,j}\,\chi _{Q_{j}}(x,t),\quad 
\kappa_{e}(x,t)=\sum\limits_{j=1}^{m}u^*_{e,j}\,\chi _{Q_{j}}(x,t), \non\\
&\quad \kappa_{i}(x,t) =\sum\limits_{j=1}^{m}u^*_{i,j}\,\chi _{Q_{j}}(x,t),\quad
\kappa_{r}(x,t)=\sum\limits_{j=1}^{m}u^*_{r,j}\, \chi _{Q_{j}}(x,t), \quad (x,t)\in Q \label{pier4}
\end{align}%
and let $\soluzstar$ denote the corresponding state. Also, let $\soluza$ be the solution to the adjoint problem \Pbla. Then, \gagi{noting that $\intQ \chi_{Q_j}= |\Omega_j|\,T$}, setting 
\begin{align}
&\mu^s_j := \frac {\int_Q (\nabla\sstar \cdot \nabla p) \, \chi_{Q_j}} {\gagi{|\Omega_j|\,T}}  , \quad 
\mu^e_j := \frac {\int_Q (\nabla\estar \cdot \nabla q )\, \chi_{Q_j}} {\gagi{|\Omega_j|\,T}}  ,\non \\[2mm]
&\quad \mu^i_j := \frac {\int_Q (\nabla\istar \cdot \nabla w )\, \chi_{Q_j}} {\gagi{|\Omega_j|\,T}} , \quad 
\mu^r_j := \frac {\int_Q (\nabla\rstar \cdot \nabla z )\, \chi_{Q_j}} {\gagi{|\Omega_j|\,T}}
 \label{pier5}
\end {align}
and recalling that $[\umin^{s,j},\umax^{s,j}]$, $[\umin^{e,j},\umax^{e,j}]$, $[\umin^{i,j},\umax^{i,j}]$, $[\umin^{r,j},\umax^{r,j}]$ are the variability intervals, all contained in $[\kamin,\kamax]$, \gagi{for $j=1,\dots,m$,} it turns out that}
\begin{align}
  & \hbox{$u^*_{s,j}$ is the one-dimensional projection of $\mu^s_j/\alpha$ on $[\umin^{s,j},\umax^{s,j}]$, i.e.,}
  \non
  \\
  & \quad u^*_{s,j} = \max\graffe{\min\graffe{\umax^{s,j},\mu^s_j/\alpha},\umin^{s,j}}
  \quad \hbox{for $j=1,\dots,m$}
  \label{pier6}
  \\[2mm]
  & \hbox{$u^*_{e,j}$ is the one-dimensional projection of $\mu^e_j/\alpha$ on $[\umin^{e,j},\umax^{e,j}]$, i.e.,}
  \non
  \\
  & \quad u^*_{e,j} = \max\graffe{\min\graffe{\umax^{e,j},\mu^e_j/\alpha},\umin^{e,j}}
  \quad \hbox{for $j=1,\dots,m$} 
  \label{pier7}
  \\[2mm]
    \separa
  & \hbox{$u^*_{i,j}$ is the one-dimensional projection of $\mu^i_j/\alpha$ on $[\umin^{i,j},\umax^{i,j}]$, i.e.,}
  \non
  \\
  &\quad u^*_{i,j} = \max\graffe{\min\graffe{\umax^{i,j},\mu^i_j/\alpha},\umin^{i,j}}
  \quad \hbox{for $j=1,\dots,m$} 
  \label{pier8}
  \\[2mm]
  & \hbox{$u^*_{r,j}$ is the one-dimensional projection of $\mu^r_j/\alpha$ on $[\umin^{r,j},\umax^{r,j}]$, i.e.,}
  \non
  \\
  & \quad u^*_{,j} = \max\graffe{\min\graffe{\umax^{r,j},\mu^r_j/\alpha},\umin^{r,j}}
  \quad \hbox{for $j=1,\dots,m$.} 
  \label{pier9}
\end{align}
\Ethm
\Bdim
We fix $\control\in\Uad$ and consider the associated linearized system \Pbll.
We test the equations by $p$, $q$, $w$ and~$z$, respectively,
and integrate with respect to time over~$(0,T)$.
At the same time, we test \accorpa{primaa}{quartaa} by $-\xi$, $-\eta$, $-\iota$ and~$-\rho$, respectively,
and integrate over~$(0,T)$.
Finally, we add all the resulting equality to each other and rearrange.
Due to the obvious cancellations that occur, we have~that
\begin{align}
  & \ioT ( \< \dt\xi(t) , p(t) > + \< \dt p(t) , \xi(t) > ) \, dt
  + \ioT ( \< \dt\eta(t) , q(t) > + \< \dt q(t) , \eta(t) > ) \, dt
  \non
  \\
  & {}
  + \ioT ( \< \dt\iota(t) , w(t) > + \< \dt w(t) , \iota(t) > ) \, dt
  + \ioT ( \< \dt\rho(t) , z(t) > + \< \dt z(t) , \rho(t) > ) \, dt
  \non
  \\
  &  {}
  \pier{={}} - \intQ (\kas-\kasstar) \nabla\sstar \cdot \nabla p
  - \intQ (\kae-\kaestar) \nabla\estar \cdot \nabla q
  \non
  \\
  & \quad {}
  - \intQ (\kai-\kaistar) \nabla\istar \cdot \nabla w
  - \intQ (\kar-\karstar) \nabla\rstar \cdot \nabla z
  \non
  \\
  & \quad {}
  \pier{{}-{}} \intQC \estar \eta
  \pier{{}-{}} \intQC \istar \iota \,.
  \non
\end{align}
\pier{Note that} all the involved functions belong to $\H1\Vp\cap\pier{{}\L2V}$\pier{. Then,
owing} to the well-known integration-by-parts formula
and accounting for the initial and final conditions \eqref{cauchyl} and~\eqref{cauchya},
we deduce that the contribution due to the first two lines of the above equality vanishes.
By combining what remains with \eqref{badNC}, we obtain\pier{%
\begin{align}
  & \intQ \bigl\{ 
    (\alpha \kasstar - \nabla\sstar \cdot \nabla p) (\kas-\kasstar) 
    + (\alpha \kaestar - \nabla\estar \cdot \nabla q) (\kae-\kaestar)
  \non 
  \\
  & \quad {}
    + (\alpha \kaistar - \nabla\istar \cdot \nabla w) (\kai-\kaistar) 
    + (\alpha \karstar - \nabla\rstar \cdot \nabla z) (\kar-\karstar) 
  \bigr\} \geq 0.
  \label{goodNC}
\end{align}
Next, by \eqref{pier2} observe that $\Uad$ is the product
$\Uad=\Uad^s\times\Uad^e\times\Uad^i\times\Uad^r$
and that a control is admissible if and only if its components belong to the corresponding factors.
Then, the variables $\kas$, $\kae$, $\kai$ and $\kar$ that enter \eqref{goodNC} are as in \eqref{pier3} 
and independent from each other, so that the inequality \eqref{goodNC} splits into four independent inequalities.
In particular, we have
\Beq
  \intQ (\alpha \kasstar - \nabla\sstar \cdot \nabla p) (\kas-\kasstar) 
  \geq 0
  \quad \hbox{for every $\kas\in\Uad^s$} 
  \non
\Eeq
and, by virtue of \eqref{pier4} and \eqref{pier5}, the above inequality reduces to 
\begin{align}
  (\alpha u_{s,j}^* - \mu^s_j) (u^s_j-u_{s,j}^*) 
  \geq 0
  \quad \hbox{ for every $u^s_j\in [\umin^{s,j},\umax^{s,j}]$, for $j=1,\dots,m$.}
  \non
\end{align}
This easily leads to \eqref{pier6}. The same argumentation can be repeated for the deduction of \eqref{pier7}--\eqref{pier9}.}
\Edim

%%%%%%%%%%%%%%%%%%%%%%%%%%%%%%%%%%%%%%%%%%%%%%%%%%%%%%%%%%%%%%%%%%%%%%%%

\section*{Acknowledgments}
This research activity has been performed in the framework of the Italian-Romanian
collaboration agreement \textquotedblleft{Analysis and control of mathematical models 
for the evolution of epidemics, tumors and phase field processes}\textquotedblright\ between the
Italian CNR and the Romanian Academy.
In addition, \pier{PC gratefully mentions} some other support from the MIUR-PRIN Grant 2020F3NCPX 
``Mathematics for industry 4.0 (Math4I4)'' and the GNAMPA (Gruppo Nazionale per l'Analisi Matematica, 
la Probabilit\`a e le loro Applicazioni) of INdAM (Isti\-tuto 
Nazionale di Alta Matematica). GM acknowledges the support of a grant of the Ministry of Research,
Innovation and Digitization, CNCS - UEFISCDI, project number
\gagi{PN-III-P4-PCE-2021-0921}, within PNCDI III.

%%%%%%%%%%%%%%%%%%%%%%%%%%%%%%%%%
%% bibliography
%%%%%%%%%%%%%%%%%%%%%%%%%%%%%%%%%

%\vspace{3truemm}

\Begin{thebibliography}{10}

\pier{\bibitem{Auricchio-23} F. Auricchio, P. Colli, G. Gilardi, A. Reali, E. Rocca,
Well-posedness for a diffusion-reaction compartmental model simulating the spread of COVID-19,
{{\it  Math. Methods Appl. Sci.} {\bf 46} (2023) 12529--12548.}}

\pier{\bibitem{Berestycki} H. Berestycki, B. Desjardins, J.S. Weitz, J.M. Oury,
Epidemic modeling with heterogeneity and social diffusion, {\it J. Math. Biol.}
{\bf 86} (2023) Paper No. 60, 59 pp.}

\bibitem{CGMR5}
P. Colli, G. Gilardi, G. Marinoschi, E. Rocca,
Optimal control of a reaction-diffusion model related to the spread of COVID-19,
\gagi{{\it Anal. Appl. (Singap.)}}, to appear  (see also preprint arXiv:2304.11114 [math.OC] (2023) pp.~1-24).

\pier{%
\bibitem{Mottoni-Orlandi-Tesei} P. de Mottoni, E. Orlandi, A. Tesei,
Asymptotic behavior for a system describing epidemics with migration and
spatial spread of infection, {\it Nonlinear Anal.}  {\bf 3}  (1979) 663-675.}

\pier{\bibitem{dOnofrio} A. d'Onofrio, M. Iannelli, P. Manfredi, G. Marinoschi,
Optimal epidemic control by social distancing and vaccination of an
infection structured by time since infection: The COVID-19 case study, {\it SIAM
J. Appl. Math.} S199-S224. https://doi.org/10.1137/22M1499406.}

\revis{%
\bibitem{morgan} W.E. Fitzgibbon, J.J. Morgan, B.Q. Tang, H.-M. Yin,
Reaction-diffusion-advection systems with discontinuous diffusion and mass control,
 {\it SIAM J. Math. Anal.} {\bf 53} (2021) 6771-6803.}

\pier{%
\bibitem{Fitzgibbon-2} W.E. Fitzgibbon, J.J. Morgan, G.F. Webb, Y. Wu, 
A vector-host epidemic model with spatial structure and age of infection,
{\it Nonlinear Anal. Real World Appl.} {\bf 41} (2018) 692-705.}

\pier{%
\bibitem{Fitzgibbon-1} W.E. Fitzgibbon, M.E. Parrott, G.F. Webb, Diffusion
epidemic models with incubation and crisscross dynamics, {\it Math.
Biosci.} {\bf 128} (1995) 131-155.}

\pier{%
\bibitem{Giordano}
G. Giordano, F. Blanchini, R. Bruno, P. Colaneri, A. Di Filippo, A. Di Matteo, M. Colaneri,
Modelling the COVID-19 epidemic and implementation of population-wide interventions in Italy,
{\it Nat. Med.} {\bf 26} (2020) 855-860.}

%\bibitem{LSU}
%O.\,A. Lady\v{z}enskaja, V.\,A. \pier{Solonnikov, N.\,N.} Uralceva,
%``Linear and Quasilinear Equations of Parabolic
%Type'', {\it Mathematical Monographs}, Vol.~{\bf 23}, 
%American Mathematical Society, Providence, Rhode Island, 1968.

\pier{\bibitem{LLLO}
W. Lee, 
S. Liu, 
W. Li, 
S. Osher, 
Mean field control problems for vaccine distribution,
{\it Res. Math. Sci.}
{\bf 9}
(2022)
Paper No. 51, 33 pp.}

\pier{\bibitem{LLTLO}
W. Lee, 
S. Liu, 
H. Tembine,
W. Li, 
S. Osher, 
Controlling propagation of epidemics via mean-field control,
{\it SIAM J. Appl. Math.}
{\bf 81}
(2021)
190-207.}

\bibitem{Lions}
J.-L. Lions, ``Quelques M\'ethodes de R\'esolution des Probl\`emes aux Limites non Lin\'eaires'', 
\pier{Dunod, Paris; Gauthier-Villars, Paris}, 1969. 

\pier{\bibitem{ManOno} P. Manfredi and A. d'Onofrio eds, 
``Modeling the interplay between human behavior and the spread of infectious diseases'',
Springer, New York, 2013.}

\pier{%
\bibitem{GM-AMO} G. Marinoschi, Parameter estimation of an epidemic model
with state constraints, \pier{{\it Appl. Math. Optim.} {\bf 84} (2021) suppl. 2, S1903-S1923.}}

\pier{%
\bibitem{GM-DCDS} G. Marinoschi, Identification of transmission rates and
reproduction number in a SARS-CoV-2 epidemic model, {\it Discrete Contin. Dyn.
Syst. Ser. S} {\bf 15} (2022) 3735-3744.}

\pier{%
\bibitem{Medhaoui} M. Mehdaoui, A.L. Alaoui, M. Tilioua,  Optimal control
for a multi-group reaction-diffusion SIR model with heterogeneous incidence
rates, {\it Int. J. Dynam. Control} {{\bf 11} (2023) 1310-1329.}}

\pier{\bibitem{Murray} J.D. Murray,  ``Mathematical biology. Second edition'',
{\it Biomathematics}, {\bf 19}. Springer-Verlag, Berlin, 1993.}

\pier{\bibitem{parino} F. Parino, L. Zino, G.C. Calafiore, A. Rizzo, A model predictive control approach to optimally devise a two-dose vaccination rollout: a case study on COVID-19 in Italy,
{\it Internat. J. Robust Nonlinear Control}  {\bf 33}  (2023) 4808-4823.}

\pier{\bibitem{Riley} S. Riley, Large-scale spatial-transmission models of
infectious disease, {\it Science}  (2007) 316(5829):1298-301.}

\pier{\bibitem{Song-Xiao} P. Song, Y. Xiao,
Analysis of a diffusive epidemic system with spatial heterogeneity and lag effect of media impact,
{\it J. Math. Biol.}  {\bf 85}  (2022) Paper No. 17, 33 pp.}

\pier{\bibitem{Viguerie-21} A. Viguerie, G. Lorenzo, F. Auricchio, D. Baroli,
T.J.R. Hughes, A. Patton, A. Reali,
Th.E. Yankeelov, A. Veneziani,
Simulating the spread of COVID-19 via a spatially-resolved
susceptible-exposed-infected-recovered-deceased (SEIRD) model with heterogeneous diffusion,
{\it Appl. Math. Lett.} {\bf 111} (2021) Paper No. 106617, 9 pp.}

\pier{\bibitem{Viguerie-20} A. Viguerie, A. Veneziani, G. Lorenzo, D. Baroli, N.
Aretz-Nellesen, A. Patton, T.E. Yankeelov, A. Reali, T.J.R. Hughes, F.
Auricchio, Diffusion-reaction compartmental models formulated in a continuum
mechanics framework: application to COVID-19, mathematical analysis, and
numerical study, {\it Comput. Mech.} {\bf 66} (2020) 1131-1152.}

\pier{\bibitem{vija} G.M. Vijayalakshmi, P. Roselyn Besi,
Vaccination control measures of an epidemic model with long-term memristive effect,
{\it J. Comput. Appl. Math.}  {\bf 419}  (2023) Paper No. 114738, 14 pp.}

\pier{%
\bibitem{Webb} G.F. Webb, A reaction-diffusion model for a deterministic
diffusive epidemic, {\it J. Math. Anal. Appl.} {\bf 84} (1981) 150-161.}

\pier{\bibitem{xu}  C. Xu, X. Huang, J. Cui, Z. Zhang, ZY. Feng, k. Cheng, 
Meta-population model about immigrants and natives with heterogeneity 
mixing and vaccine strategy of tuberculosis in China,
{\it Int. J. Biomath.} {\bf 16}  (2023) Paper No. 2250121, 10 pp.}

\End{thebibliography}

\End{document}

%%%%%%%%%%%%%%%%%%%%%%%%%%%%%%%%%%%%%%%%%%%%%